\documentclass[8pt,a4wide]{article}
\usepackage[numbers,sort&compress]{natbib}
\usepackage[normalem]{ulem}
\usepackage{listings}
\usepackage{showlabels}
\usepackage{a4wide}
\usepackage{amsmath} 
\usepackage{amsthm}
\usepackage{amsfonts}
\usepackage{amssymb,amscd}
\usepackage{multirow}
\usepackage[mathcal]{eucal}
\usepackage[spanish,USenglish]{babel} 
\usepackage{graphics}  
\usepackage{graphicx} 
\usepackage{hyperref}
\usepackage{fancyhdr}
\usepackage{titlesec}
\usepackage[latin1]{inputenc}
\usepackage{amsmath}
\usepackage{amsfonts}
\usepackage{amssymb}
\usepackage{amsthm}
\usepackage{newlfont}
\usepackage{amsbsy}
\usepackage{bm}
\usepackage[all]{xy}
\usepackage{textcomp}
\usepackage{color}
\usepackage{epsfig}
\usepackage{dsfont}
\usepackage{latexsym}
\usepackage{array}
\usepackage{longtable}
\usepackage{enumerate}
\usepackage{multicol}
\usepackage{mathrsfs}
\usepackage{cancel}
\usepackage{wasysym}
\theoremstyle{plain}
\newtheorem{teor}{Theorem}
\newtheorem{lemma}{Lemma}

\newtheorem{prop}{Proposition}

\theoremstyle{definition}

\newtheorem*{ejem*}{Examples}
\newtheorem{defin}{Definition}
\newtheorem*{demos}{Proof}
\newtheorem{remark}{Remark}
\theoremstyle{remark}

\newcommand{\complex}{\mathbf{\mathbb{C}}}

\newcommand{\natu}{\mathbf{\mathbb{N}}}
\newcommand{\llave}[1]{\left\{ #1\right\}}

\newcommand{\corch}[1]{\left[ #1\right]}
\newcommand{\paren}[1]{\left( #1\right)}

\usepackage[breakable]{tcolorbox}
\usepackage{parskip} 

\usepackage{iftex}
\ifPDFTeX
\usepackage[T1]{fontenc}
\usepackage{mathpazo}
\else
\usepackage{fontspec}
\fi

\usepackage{graphicx}

\usepackage{caption}
\DeclareCaptionFormat{nocaption}{}
\usepackage[Export]{adjustbox} 
\adjustboxset{max size={0.9\linewidth}{0.9\paperheight}}
\usepackage{float}
\usepackage{xcolor} 
\usepackage{enumerate} 
\usepackage{geometry} 
\usepackage{amsmath} 
\usepackage{amssymb} 
\usepackage{textcomp} 
\AtBeginDocument{%
}
\usepackage{upquote} 
\usepackage{eurosym} 
\usepackage[mathletters]{ucs} 
\usepackage{fancyvrb} 
\usepackage{grffile} 
\makeatletter 
\def\Gread@@xetex#1{%
	\IfFileExists{"\Gin@base".bb}%
	{\Gread@eps{\Gin@base.bb}}%
	{\Gread@@xetex@aux#1}%
}
\makeatother

\usepackage{hyperref}
\usepackage{titling}
\usepackage{longtable} 
\usepackage{booktabs}  
\usepackage[inline]{enumitem} 
\usepackage[normalem]{ulem} 
\usepackage{mathrsfs}

\definecolor{urlcolor}{rgb}{0,.145,.698}
\definecolor{linkcolor}{rgb}{.71,0.21,0.01}
\definecolor{citecolor}{rgb}{.12,.54,.11}

\definecolor{ansi-black}{HTML}{3E424D}
\definecolor{ansi-black-intense}{HTML}{282C36}
\definecolor{ansi-red}{HTML}{E75C58}
\definecolor{ansi-red-intense}{HTML}{B22B31}
\definecolor{ansi-green}{HTML}{00A250}
\definecolor{ansi-green-intense}{HTML}{007427}
\definecolor{ansi-yellow}{HTML}{DDB62B}
\definecolor{ansi-yellow-intense}{HTML}{B27D12}
\definecolor{ansi-blue}{HTML}{208FFB}
\definecolor{ansi-blue-intense}{HTML}{0065CA}
\definecolor{ansi-magenta}{HTML}{D160C4}
\definecolor{ansi-magenta-intense}{HTML}{A03196}
\definecolor{ansi-cyan}{HTML}{60C6C8}
\definecolor{ansi-cyan-intense}{HTML}{258F8F}
\definecolor{ansi-white}{HTML}{C5C1B4}
\definecolor{ansi-white-intense}{HTML}{A1A6B2}
\definecolor{ansi-default-inverse-fg}{HTML}{FFFFFF}
\definecolor{ansi-default-inverse-bg}{HTML}{000000}

\let\Oldtex\TeX
\let\Oldlatex\LaTeX
\renewcommand{\TeX}{\textrm{\Oldtex}}
\renewcommand{\LaTeX}{\textrm{\Oldlatex}}
\title{First Example}

\makeatletter
\def\PY@reset{\let\PY@it=\relax \let\PY@bf=\relax%
	\let\PY@ul=\relax \let\PY@tc=\relax%
	\let\PY@bc=\relax \let\PY@ff=\relax}
\def\PY@tok#1{\csname PY@tok@#1\endcsname}
\def\PY@toks#1+{\ifx\relax#1\empty\else%
	\PY@tok{#1}\expandafter\PY@toks\fi}
\def\PY@do#1{\PY@bc{\PY@tc{\PY@ul{%
				\PY@it{\PY@bf{\PY@ff{#1}}}}}}}
\def\PY#1#2{\PY@reset\PY@toks#1+\relax+\PY@do{#2}}

\expandafter\def\csname PY@tok@w\endcsname{\def\PY@tc##1{\textcolor[rgb]{0.73,0.73,0.73}{##1}}}
\expandafter\def\csname PY@tok@c\endcsname{\let\PY@it=\textit\def\PY@tc##1{\textcolor[rgb]{0.25,0.50,0.50}{##1}}}
\expandafter\def\csname PY@tok@cp\endcsname{\def\PY@tc##1{\textcolor[rgb]{0.74,0.48,0.00}{##1}}}
\expandafter\def\csname PY@tok@k\endcsname{\let\PY@bf=\textbf\def\PY@tc##1{\textcolor[rgb]{0.00,0.50,0.00}{##1}}}
\expandafter\def\csname PY@tok@kp\endcsname{\def\PY@tc##1{\textcolor[rgb]{0.00,0.50,0.00}{##1}}}
\expandafter\def\csname PY@tok@kt\endcsname{\def\PY@tc##1{\textcolor[rgb]{0.69,0.00,0.25}{##1}}}
\expandafter\def\csname PY@tok@o\endcsname{\def\PY@tc##1{\textcolor[rgb]{0.40,0.40,0.40}{##1}}}
\expandafter\def\csname PY@tok@ow\endcsname{\let\PY@bf=\textbf\def\PY@tc##1{\textcolor[rgb]{0.67,0.13,1.00}{##1}}}
\expandafter\def\csname PY@tok@nb\endcsname{\def\PY@tc##1{\textcolor[rgb]{0.00,0.50,0.00}{##1}}}
\expandafter\def\csname PY@tok@nf\endcsname{\def\PY@tc##1{\textcolor[rgb]{0.00,0.00,1.00}{##1}}}
\expandafter\def\csname PY@tok@nc\endcsname{\let\PY@bf=\textbf\def\PY@tc##1{\textcolor[rgb]{0.00,0.00,1.00}{##1}}}
\expandafter\def\csname PY@tok@nn\endcsname{\let\PY@bf=\textbf\def\PY@tc##1{\textcolor[rgb]{0.00,0.00,1.00}{##1}}}
\expandafter\def\csname PY@tok@ne\endcsname{\let\PY@bf=\textbf\def\PY@tc##1{\textcolor[rgb]{0.82,0.25,0.23}{##1}}}
\expandafter\def\csname PY@tok@nv\endcsname{\def\PY@tc##1{\textcolor[rgb]{0.10,0.09,0.49}{##1}}}
\expandafter\def\csname PY@tok@no\endcsname{\def\PY@tc##1{\textcolor[rgb]{0.53,0.00,0.00}{##1}}}
\expandafter\def\csname PY@tok@nl\endcsname{\def\PY@tc##1{\textcolor[rgb]{0.63,0.63,0.00}{##1}}}
\expandafter\def\csname PY@tok@ni\endcsname{\let\PY@bf=\textbf\def\PY@tc##1{\textcolor[rgb]{0.60,0.60,0.60}{##1}}}
\expandafter\def\csname PY@tok@na\endcsname{\def\PY@tc##1{\textcolor[rgb]{0.49,0.56,0.16}{##1}}}
\expandafter\def\csname PY@tok@nt\endcsname{\let\PY@bf=\textbf\def\PY@tc##1{\textcolor[rgb]{0.00,0.50,0.00}{##1}}}
\expandafter\def\csname PY@tok@nd\endcsname{\def\PY@tc##1{\textcolor[rgb]{0.67,0.13,1.00}{##1}}}
\expandafter\def\csname PY@tok@s\endcsname{\def\PY@tc##1{\textcolor[rgb]{0.73,0.13,0.13}{##1}}}
\expandafter\def\csname PY@tok@sd\endcsname{\let\PY@it=\textit\def\PY@tc##1{\textcolor[rgb]{0.73,0.13,0.13}{##1}}}
\expandafter\def\csname PY@tok@si\endcsname{\let\PY@bf=\textbf\def\PY@tc##1{\textcolor[rgb]{0.73,0.40,0.53}{##1}}}
\expandafter\def\csname PY@tok@se\endcsname{\let\PY@bf=\textbf\def\PY@tc##1{\textcolor[rgb]{0.73,0.40,0.13}{##1}}}
\expandafter\def\csname PY@tok@sr\endcsname{\def\PY@tc##1{\textcolor[rgb]{0.73,0.40,0.53}{##1}}}
\expandafter\def\csname PY@tok@ss\endcsname{\def\PY@tc##1{\textcolor[rgb]{0.10,0.09,0.49}{##1}}}
\expandafter\def\csname PY@tok@sx\endcsname{\def\PY@tc##1{\textcolor[rgb]{0.00,0.50,0.00}{##1}}}
\expandafter\def\csname PY@tok@m\endcsname{\def\PY@tc##1{\textcolor[rgb]{0.40,0.40,0.40}{##1}}}
\expandafter\def\csname PY@tok@gh\endcsname{\let\PY@bf=\textbf\def\PY@tc##1{\textcolor[rgb]{0.00,0.00,0.50}{##1}}}
\expandafter\def\csname PY@tok@gu\endcsname{\let\PY@bf=\textbf\def\PY@tc##1{\textcolor[rgb]{0.50,0.00,0.50}{##1}}}
\expandafter\def\csname PY@tok@gd\endcsname{\def\PY@tc##1{\textcolor[rgb]{0.63,0.00,0.00}{##1}}}
\expandafter\def\csname PY@tok@gi\endcsname{\def\PY@tc##1{\textcolor[rgb]{0.00,0.63,0.00}{##1}}}
\expandafter\def\csname PY@tok@gr\endcsname{\def\PY@tc##1{\textcolor[rgb]{1.00,0.00,0.00}{##1}}}
\expandafter\def\csname PY@tok@ge\endcsname{\let\PY@it=\textit}
\expandafter\def\csname PY@tok@gs\endcsname{\let\PY@bf=\textbf}
\expandafter\def\csname PY@tok@gp\endcsname{\let\PY@bf=\textbf\def\PY@tc##1{\textcolor[rgb]{0.00,0.00,0.50}{##1}}}
\expandafter\def\csname PY@tok@go\endcsname{\def\PY@tc##1{\textcolor[rgb]{0.53,0.53,0.53}{##1}}}
\expandafter\def\csname PY@tok@gt\endcsname{\def\PY@tc##1{\textcolor[rgb]{0.00,0.27,0.87}{##1}}}
\expandafter\def\csname PY@tok@err\endcsname{\def\PY@bc##1{\setlength{\fboxsep}{0pt}\fcolorbox[rgb]{1.00,0.00,0.00}{1,1,1}{\strut ##1}}}
\expandafter\def\csname PY@tok@kc\endcsname{\let\PY@bf=\textbf\def\PY@tc##1{\textcolor[rgb]{0.00,0.50,0.00}{##1}}}
\expandafter\def\csname PY@tok@kd\endcsname{\let\PY@bf=\textbf\def\PY@tc##1{\textcolor[rgb]{0.00,0.50,0.00}{##1}}}
\expandafter\def\csname PY@tok@kn\endcsname{\let\PY@bf=\textbf\def\PY@tc##1{\textcolor[rgb]{0.00,0.50,0.00}{##1}}}
\expandafter\def\csname PY@tok@kr\endcsname{\let\PY@bf=\textbf\def\PY@tc##1{\textcolor[rgb]{0.00,0.50,0.00}{##1}}}
\expandafter\def\csname PY@tok@bp\endcsname{\def\PY@tc##1{\textcolor[rgb]{0.00,0.50,0.00}{##1}}}
\expandafter\def\csname PY@tok@fm\endcsname{\def\PY@tc##1{\textcolor[rgb]{0.00,0.00,1.00}{##1}}}
\expandafter\def\csname PY@tok@vc\endcsname{\def\PY@tc##1{\textcolor[rgb]{0.10,0.09,0.49}{##1}}}
\expandafter\def\csname PY@tok@vg\endcsname{\def\PY@tc##1{\textcolor[rgb]{0.10,0.09,0.49}{##1}}}
\expandafter\def\csname PY@tok@vi\endcsname{\def\PY@tc##1{\textcolor[rgb]{0.10,0.09,0.49}{##1}}}
\expandafter\def\csname PY@tok@vm\endcsname{\def\PY@tc##1{\textcolor[rgb]{0.10,0.09,0.49}{##1}}}
\expandafter\def\csname PY@tok@sa\endcsname{\def\PY@tc##1{\textcolor[rgb]{0.73,0.13,0.13}{##1}}}
\expandafter\def\csname PY@tok@sb\endcsname{\def\PY@tc##1{\textcolor[rgb]{0.73,0.13,0.13}{##1}}}
\expandafter\def\csname PY@tok@sc\endcsname{\def\PY@tc##1{\textcolor[rgb]{0.73,0.13,0.13}{##1}}}
\expandafter\def\csname PY@tok@dl\endcsname{\def\PY@tc##1{\textcolor[rgb]{0.73,0.13,0.13}{##1}}}
\expandafter\def\csname PY@tok@s2\endcsname{\def\PY@tc##1{\textcolor[rgb]{0.73,0.13,0.13}{##1}}}
\expandafter\def\csname PY@tok@sh\endcsname{\def\PY@tc##1{\textcolor[rgb]{0.73,0.13,0.13}{##1}}}
\expandafter\def\csname PY@tok@s1\endcsname{\def\PY@tc##1{\textcolor[rgb]{0.73,0.13,0.13}{##1}}}
\expandafter\def\csname PY@tok@mb\endcsname{\def\PY@tc##1{\textcolor[rgb]{0.40,0.40,0.40}{##1}}}
\expandafter\def\csname PY@tok@mf\endcsname{\def\PY@tc##1{\textcolor[rgb]{0.40,0.40,0.40}{##1}}}
\expandafter\def\csname PY@tok@mh\endcsname{\def\PY@tc##1{\textcolor[rgb]{0.40,0.40,0.40}{##1}}}
\expandafter\def\csname PY@tok@mi\endcsname{\def\PY@tc##1{\textcolor[rgb]{0.40,0.40,0.40}{##1}}}
\expandafter\def\csname PY@tok@il\endcsname{\def\PY@tc##1{\textcolor[rgb]{0.40,0.40,0.40}{##1}}}
\expandafter\def\csname PY@tok@mo\endcsname{\def\PY@tc##1{\textcolor[rgb]{0.40,0.40,0.40}{##1}}}
\expandafter\def\csname PY@tok@ch\endcsname{\let\PY@it=\textit\def\PY@tc##1{\textcolor[rgb]{0.25,0.50,0.50}{##1}}}
\expandafter\def\csname PY@tok@cm\endcsname{\let\PY@it=\textit\def\PY@tc##1{\textcolor[rgb]{0.25,0.50,0.50}{##1}}}
\expandafter\def\csname PY@tok@cpf\endcsname{\let\PY@it=\textit\def\PY@tc##1{\textcolor[rgb]{0.25,0.50,0.50}{##1}}}
\expandafter\def\csname PY@tok@c1\endcsname{\let\PY@it=\textit\def\PY@tc##1{\textcolor[rgb]{0.25,0.50,0.50}{##1}}}
\expandafter\def\csname PY@tok@cs\endcsname{\let\PY@it=\textit\def\PY@tc##1{\textcolor[rgb]{0.25,0.50,0.50}{##1}}}

\makeatother

\makeatletter
\newbox\Wrappedcontinuationbox 
\newbox\Wrappedvisiblespacebox 
\newcommand*\Wrappedvisiblespace {\textcolor{red}{\textvisiblespace}} 
\newcommand*\Wrappedcontinuationsymbol {\textcolor{red}{\llap{\tiny$\m@th\hookrightarrow$}}} 
\newcommand*\Wrappedcontinuationindent {3ex } 
\newcommand*\Wrappedafterbreak {\kern\Wrappedcontinuationindent\copy\Wrappedcontinuationbox} 
\newcommand*\Wrappedbreaksatspecials {%
	\def\PYGZus{\discretionary{\char`\_}{\Wrappedafterbreak}{\char`\_}}%
	\def\PYGZob{\discretionary{}{\Wrappedafterbreak\char`\{}{\char`\{}}%
	\def\PYGZcb{\discretionary{\char`\}}{\Wrappedafterbreak}{\char`\}}}%
	\def\PYGZca{\discretionary{\char`\^}{\Wrappedafterbreak}{\char`\^}}%
	\def\PYGZam{\discretionary{\char`\&}{\Wrappedafterbreak}{\char`\&}}%
	\def\PYGZlt{\discretionary{}{\Wrappedafterbreak\char`\<}{\char`\<}}%
	\def\PYGZgt{\discretionary{\char`\>}{\Wrappedafterbreak}{\char`\>}}%
	\def\PYGZsh{\discretionary{}{\Wrappedafterbreak\char`\#}{\char`\#}}%
	\def\PYGZpc{\discretionary{}{\Wrappedafterbreak\char`\%}{\char`\%}}%
	\def\PYGZdl{\discretionary{}{\Wrappedafterbreak\char`\$}{\char`\$}}%
	\def\PYGZhy{\discretionary{\char`\-}{\Wrappedafterbreak}{\char`\-}}%
	\def\PYGZsq{\discretionary{}{\Wrappedafterbreak\textquotesingle}{\textquotesingle}}%
	\def\PYGZdq{\discretionary{}{\Wrappedafterbreak\char`\"}{\char`\"}}%
	\def\PYGZti{\discretionary{\char`\~}{\Wrappedafterbreak}{\char`\~}}%
} 
\newcommand*\Wrappedbreaksatpunct {%
	\lccode`\~`\.\lowercase{\def~}{\discretionary{\hbox{\char`\.}}{\Wrappedafterbreak}{\hbox{\char`\.}}}%
	\lccode`\~`\,\lowercase{\def~}{\discretionary{\hbox{\char`\,}}{\Wrappedafterbreak}{\hbox{\char`\,}}}%
	\lccode`\~`\;\lowercase{\def~}{\discretionary{\hbox{\char`\;}}{\Wrappedafterbreak}{\hbox{\char`\;}}}%
	\lccode`\~`\:\lowercase{\def~}{\discretionary{\hbox{\char`\:}}{\Wrappedafterbreak}{\hbox{\char`\:}}}%
	\lccode`\~`\?\lowercase{\def~}{\discretionary{\hbox{\char`\?}}{\Wrappedafterbreak}{\hbox{\char`\?}}}%
	\lccode`\~`\!\lowercase{\def~}{\discretionary{\hbox{\char`\!}}{\Wrappedafterbreak}{\hbox{\char`\!}}}%
	\lccode`\~`\/\lowercase{\def~}{\discretionary{\hbox{\char`\/}}{\Wrappedafterbreak}{\hbox{\char`\/}}}%
	\catcode`\.\active
	\catcode`\,\active 
	\catcode`\;\active
	\catcode`\:\active
	\catcode`\?\active
	\catcode`\!\active
	\catcode`\/\active 
	\lccode`\~`\~ 	
}
\makeatother

\let\OriginalVerbatim=\Verbatim
\makeatletter
\renewcommand{\Verbatim}[1][1]{%
	\sbox\Wrappedcontinuationbox {\Wrappedcontinuationsymbol}%
	\sbox\Wrappedvisiblespacebox {\FV@SetupFont\Wrappedvisiblespace}%
	\def\FancyVerbFormatLine ##1{\hsize\linewidth
		\vtop{\raggedright\hyphenpenalty\z@\exhyphenpenalty\z@
			\doublehyphendemerits\z@\finalhyphendemerits\z@
			\strut ##1\strut}%
	}%
	\def\FV@Space {%
		\nobreak\hskip\z@ plus\fontdimen3\font minus\fontdimen4\font
		\discretionary{\copy\Wrappedvisiblespacebox}{\Wrappedafterbreak}
		{\kern\fontdimen2\font}%
	}%
	
	\Wrappedbreaksatspecials
	\OriginalVerbatim[#1,codes*=\Wrappedbreaksatpunct]%
}
\makeatother

\definecolor{incolor}{HTML}{303F9F}
\definecolor{outcolor}{HTML}{D84315}
\definecolor{cellborder}{HTML}{CFCFCF}
\definecolor{cellbackground}{HTML}{F7F7F7}

\makeatletter
\newcommand{\boxspacing}{\kern\kvtcb@left@rule\kern\kvtcb@boxsep}
\makeatother

\hypersetup{
	breaklinks=true,  
	colorlinks=true,
	urlcolor=urlcolor,
	linkcolor=linkcolor,
	citecolor=citecolor,
}

\begin{document}
	\title{\textbf{Noncommutative Bispectral Algebras and their Presentations}}
	\author{Brian D. Vasquez C. & Jorge P. Zubelli}
	\author{
Brian D. Vasquez$^{1}$, Jorge P. Zubelli$^{2}$\\ \\
\small{$^{1}$IMPA, $^{2}$Department of Mathematics, Khalifa University.}\\
\small{\texttt{$^{1}$bridava927@gmail.com, $^{2}$zubelli@gmail.com}}
}
\date{\small{\today}}
	\maketitle

\begin{abstract}
We prove a general result on presentations of finitely-generated algebras and apply it to obtain nice presentations for some noncommutative  algebras arising in the matrix bispectral problem.  By ``nice presentation'' we mean
a presentation that has as few as possible defining relations. This in turn, has potential applications in computer algebra implementations and examples. 

Our results can be divided into three parts. 
 In the first two, we consider bispectral algebras with the eigenvalue in the physical equation to be scalar-valued for $2\times 2$ and $3\times 3$ matrix-valued eigenfunctions. In the third part,  we assume the eigenvalue in the physical equation to be matrix-valued and draw an important connection with spin Calogero-Moser systems.
 In all cases, we show that these algebras are finitely presented.
As a byproduct, we answer positively  a  conjecture of F.~A.~Gr{\"{u}}nbaum about these algebras.\\
\end{abstract}

\emph{Key words}: bispectral problem, Calogero-Moser systems, presentations of finitely generated algebras, completely integrable systems.

\linespread{2}

\section{Introduction}
In this work, we characterize the symmetry structure of a noncommutative version of the bispectral problem \cite{duistermaat1986differential}.
The latter refers to families of eigenfunctions
$\psi (x,z)$ of an operator $L=L(x,\partial_{x})$, with  $z$-dependent eigenvalue parameter, that are also eigenfunctions for some nontrivial operator $B=B(z,\partial_{z})$ with an  $x$-dependent eigenvalue. 

In the scalar case, the bispectral problem already displays unexpected connections to different areas \cite{duistermaat1986differential,Zubelli1991,marta2018}. One of the most important connections is that a remarkable set of bispectral Schr\"{o}dinger operators $L=-\partial_{x}^2+U(x)$ are obtained when $U(x)$ is a rational solution of the KdV equation \cite{AMM1977}. The abundance of connections is even more pronounced in the matrix case. See \cite{Zubelli1991,wilson,Chalub2000,Chalub2001,Chalub2001a,MR2201201,Sakhnovich2001,PNAS2019} and references therein. In the theory of infinite dimensional systems and solitons the study of the symmetries led to a deeper  understanding of the structure of these equations. See for example~\cite{fokas1987,fokas2002}.

Characterizing the algebraic structure of the solutions to a problem through presentations is a major task in many areas. In our context, 
this consists in looking for a set of generators in such a way that the relations among them are as simple as possible \cite{DERKSEN2015210,zbMATH03512288,zbMATH01713734}.  
We address this problem for some algebras associated to the noncommutative bispectral problem~\cite{duistermaat1986differential}.

Interesting conjectures concerning presentations of some noncommutative algebras were proposed in connection with the interplay of   matrix-valued orthogonal polynomials \cite{Grfrm-4,zbMATH05249007} and the bispectral problem~\cite{zbMATH01684157}.
Only one of the conjectures proposed in \cite{zbMATH05249007} was solved in \cite{tirao}. In \cite{Grfrm-4} the algebras involved are bispectral algebras while in \cite{zbMATH05249007} the algebras involved are algebras of differential operators associated to matrix-valued orthogonal polynomials. This article solves the conjectures concerning noncommutative bispectral algebras presented in \cite{Grfrm-4}.

In the present incarnation of the {\em bispectral problem}, we consider the triples $(L,\psi, B)$ satisfying systems of equations 
	\begin{equation}\label{bispec}
\left\lbrace\begin{array}{c} 
L\psi(x,z)=\psi (x,z)F(z) \\
(\psi B) (x,z)=\theta(x)\psi(x,z)
\end{array}\right.
	\end{equation}
	with $L=L(x,\partial_{x})$, $B=B(z,\partial_{z})$  linear matrix  differential operators, i.e., $L\psi=\sum_{i=0}^{l}a_{i}(x) \cdot \partial_{x}^i \psi$,
	$\psi B=\sum_{j=0}^{m}\partial_{z}^j \psi \cdot b_{j}(z)$. The functions $a_{i}, b_{j},F, \theta$ and the nontrivial common eigenfunction $\psi$ are in principle compatible sized matrix valued functions.
A triple $(L,\psi, B)$ satisfying \eqref{bispec} is called a bispectral triple. 


The main goal of this article is to give a presentation of each (bispectral) algebra using its generators and some relations among them. 
Thus, describing the ideal of relations, we give three examples of bispectral algebras to illustrate a general theorem of presentations of finitely generated algebras. For a given eigenvalue function the corresponding algebra of  matrix eigenvalues is characterized. 
In the former two cases, the  eigenvalue $F(z)$ is scalar valued and in the last case the eigenvalue $\theta(x)$ is matrix valued. These results give positive answers to the three conjectures 
in~\cite{Grfrm-4}. We use the software Singular and Maxima to obtain a set of generators and nice relations among them and after that, we prove that in fact, this set of nice relations are enough to give presentations for these algebras.

Now we fix the normalized~\footnote{ If $L=L(x,\partial_{x})$, $L=\sum_{i=0}^{l}a_{i}(x)\partial_{x}^i$ with $a_{l}$ constant and scalar, $a_{l-1}=0$, then $L$ is called normalized.} operator $L$ and the eigenfunctions $\psi(\cdot,z)$. We are interested in the bispectral pairs associated to $L=L(x,\partial_{x})$, i.e., the algebra 
\begin{equation}
\mathbb{A}=\llave{\theta \in  M_{N}(\complex)\corch{x} \big| \exists B=B(z,\partial_{z}), (\psi B)(x,z)=\theta(x)\psi(x,z)}. 
\end{equation}


In order to characterize the algebraic structure of bispectrality in 
the present noncommutative context, we start with the following definitions.
\begin{defin}
Let  $\mathbb{K}$ be a field, $C$ be a $\mathbb{K}$-algebra, $A$ a subring of $C$ and $S\subset C$. We define 
$$A\cdot <S>=span_{\mathbb{K}}\llave{\prod_{j=1}^{n}s_j\mid s_1, ... , s_n \in S\cup A, n \in \natu } \mbox{ ,}$$
where the noncommutative product is understood from left to right, i.e., 
$ \prod_{j=1}^{n+1}s_j := (\prod_{j=1}^{n}s_j )s_{n+1}, $ for $n=0,1,2,\cdots.$ For completion,
$\prod_{j=1}^{0}s_j:=1.$
\end{defin} 

The set $A\cdot<S>$ is called the {\em subalgebra generated by} $S$ over $A$ and we call an element $f\in A\cdot<S>$ a noncommutative polynomial with 
coefficients in $A$ and set of variables $S$. 

\begin{defin}
    Let $C$ be a noncommutative ring and $A$ a subring of $C$. We say that an element $\alpha \in C$ is integral over $A$ if there exists a 
    noncommutative polynomial $f$ with coefficients in $A$ such that $f(\alpha)=0$. Furthermore, we say that $\beta\in C $ is integral 
    over $\alpha \in C$ if $\beta$ is integral over $A\cdot<\alpha>$. Finally, $\alpha$ and $\beta$ are associated integrals if  $\alpha$ is integral over $\beta$ and $\beta$ is integral over $\alpha$.
\end{defin}

\begin{defin}
	The shift operator $S_{N}\in M_{N}(\mathbb{K}[x])$ is defined by  $$S_{N}=\sum_{s=1}^{N-1}e_{s,s+1}$$ 
	for $N\geq 2$, where as usual
	$e_{r,s}$ denotes the matrix with $1$ at entry $(r,s)$ and zeros elsewhere.  
\end{defin}

We consider a nilpotent element $S\in M_{N}(\mathbb{K})$ of degree $D\geq 2$,  the matrix valued function 
\begin{equation*}
\psi(x,z)=
e^{xz}\paren{Iz+\sum_{m=1}^{D}(-1)^{m} S^{m-1}x^{-m}} \mbox{ , }
\end{equation*}
and note that $L\psi(x,z)=-z^2\psi(x,z)$ for the ordinary differential operator
\begin{equation*}
L= -\partial_{x}^2+2
\sum_{m=1}^{D}(-1)^{m+1} mS^{m-1}x^{-m-1}.
\end{equation*}

In \cite{VZ1} we studied the bispectral algebra $\mathbb{A}$ associated to this operator $L$ and eigenfunction $\psi$, we give an explicit expression for the operator $B=B(z,\partial_{z})$ associated to the matrix eigenvalue $\theta$. However, we did not give a characterization in terms of generators and relations for $\mathbb{A}$. The main goal of the present work is to prove characterization results about the bispectral algebra $\mathbb{A}$ (in sizes of matrix $N=2,3$) and give positive answers to the conjectures purposed by F.~A.~Gr\"{u}nbaum in \cite{Grfrm-4}.




The following two theorems are for matrices of size $N=2,3$ and the nilpotent element is the shift operator $S=S_{N}.$


\begin{teor}[An algebra with an integral element over a nilpotent one]\label{first1}
	Let $\Gamma$ be the sub-algebra of $M_{2}(\complex)[x]$ of the form 
	\begin{equation*}
	\paren{
		\begin{matrix} 
		r_{0}^{11} & r_{0}^{12} \\
		0 & r_{0}^{11} 
		\end{matrix}}+  \paren{
		\begin{matrix} 
		r_{1}^{11} & r_{1}^{12} \\
		0 & r_{1}^{11} 
		\end{matrix}}x+ \paren{
		\begin{matrix} 
		r_{2}^{11} & r_{2}^{12} \\
		r_{1}^{11} & r_{2}^{22} 
		\end{matrix}}x^{2}+
	\paren{
		\begin{matrix} 
		r_{3}^{11} & r_{3}^{12} \\
		r_{2}^{22}+r_{2}^{11}-r_{1}^{12} & r_{3}^{22} 
		\end{matrix}}x^3+x^{4}p(x),
	\end{equation*} 
	where $p\in M_{2}(\complex)[x]$ and all the variables $r_{0}^{11},r_{0}^{12}, r_{1}^{11},r_{1}^{12},r_{2}^{11},r_{2}^{22},r_{3}^{11},r_{3}^{12}, r_{3}^{22} \in \complex$. Then $\Gamma=\mathbb{A}$. Moreover, for each $\theta$ we have an explicit expression for the operator $B$.
	
	Furthermore, we have the presentation $\mathbb{A}=\complex\cdot \langle \alpha_0, \alpha_1 \mid I=0\rangle$ with
the ideal $I$ given by 
\begin{equation*}
I :=\langle \alpha_{0}^2,\alpha_{1}^3+\alpha_{0}\alpha_{1}\alpha_{0}-3\alpha_{1}\alpha_{0}\alpha_{1}+\alpha_{0}\alpha_{1}^2+\alpha_{1}^2\alpha_{0} \rangle \mbox{ .}
\end{equation*}
\end{teor}


\begin{teor}[An algebra with nilpotent and idempotent associated elements]\label{second}

	Let $\Gamma$ the sub-algebra of $M_{3}(\complex)[x]$ of the form 
	\begin{equation*}
	\paren{
		\begin{matrix} 
		r_{0}^{11} & r_{0}^{12} & r_{0}^{13} \\
		0 & r_{0}^{22} & r_{0}^{23} \\
		0 & 0 & r_{0}^{11}
		\end{matrix}}+ \paren{
		\begin{matrix} 
		r_{1}^{11} & r_{1}^{12} & r_{1}^{13} \\
		r_{0}^{22}-r_{0}^{11} & r_{1}^{22} & r_{1}^{23} \\
		0 & r_{0}^{22}-r_{0}^{11} & r_{1}^{11}+r_{0}^{23}-r_{0}^{12}
		\end{matrix}}x
	\end{equation*}
	\begin{equation*}
	+ \paren{
		\begin{matrix} 
		r_{2}^{11} & r_{2}^{12} & r_{2}^{13} \\
		r_{1}^{22}-r_{1}^{11}-r_{0}^{23}+r_{0}^{12} & r_{2}^{22} & r_{2}^{23} \\
		r_{0}^{22}-r_{0}^{11} & r_{1}^{22}-r_{1}^{11} & r_{2}^{11}+r_{1}^{23}-r_{1}^{12}
		\end{matrix}}x^{2}
	+  \paren{
		\begin{matrix} 
		r_{3}^{11} & r_{3}^{12} & r_{3}^{13} \\
		r_{3}^{21}& r_{3}^{22} & r_{3}^{23} \\
		r_{1}^{22}-2r_{1}^{11}-r_{0}^{23}+r_{0}^{12} & r_{3}^{32} & r_{3}^{33}
		\end{matrix}}x^{3}   
	\end{equation*}
	\begin{equation*}
	+ \paren{
		\begin{matrix} 
		r_{4}^{11} & r_{4}^{12} & r_{4}^{13} \\
		r_{4}^{21} & r_{4}^{22} & r_{4}^{23} \\
		r_{3}^{32}+r_{3}^{21}-r_{2}^{22}-r_{2}^{11}+r_{1}^{12} & r_{4}^{22} & r_{4}^{33}
		\end{matrix}}x^{4}
	\end{equation*}	
	\begin{equation*}
+  \paren{
	\begin{matrix} 
	r_{5}^{11} & r_{5}^{12} & r_{5}^{13} \\
	r_{5}^{21}& r_{5}^{22} & r_{5}^{23} \\
	r_{4}^{32}+r_{4}^{21}-r_{3}^{33}-r_{3}^{22}-r_{3}^{11}+r_{2}^{23}+r_{2}^{12}-r_{1}^{13} & r_{5}^{32} & r_{5}^{33}
	\end{matrix}}x^{5}	+x^{6}p(x) \mbox{ , }
	\end{equation*} 
	where $p\in M_{3}(\complex)[x]$ and all the variables $r_{0}^{11}, r_{0}^{12},...,r_{5}^{33}\in \complex$ are arbitrary.
	
	Then, $\Gamma =\mathbb{A}$ and for each $\theta$ we have  an explicit expression for the operator $B$.
	
	Furthermore, we have the presentation $\mathbb{A}=\complex\cdot \langle \alpha_2, \alpha_3 \mid I=0\rangle$ with
	
	\begin{equation*}
	I=\langle\alpha_{2}^3,\alpha_{3}^2-\alpha_{3}, (\alpha_{3}\alpha_{2})^2 \alpha_{3}-4\alpha_{3}\alpha_{2}^2 \alpha_{3}\rangle \mbox{ .}
	\end{equation*}
\end{teor}


We now consider the case when both ``eigenvalues" $F$ and $\theta$ are matrix valued. 
Let 
\begin{equation*}
\psi(x,z)=\frac{e^{xz}}{(x-2)xz}\paren{\begin{matrix} 
\frac{x^3z^2-2x^2z^2-2x^2z+3xz+2x-2}{xz} & \frac{1}{x} \\
\frac{xz-2}{z} & x^2z-2xz-x+1
\end{matrix}}
\end{equation*}
and \begin{equation*}
L=\paren{\begin{matrix} 
     0 & 0 \\
	0 & 1
	\end{matrix}}.\partial_{x}^2
+\paren{\begin{matrix} 
0 & \frac{1}{(x-2)x^2} \\
	-\frac{1}{x-2} & 0
	\end{matrix}}.\partial_{x}
+ \paren{\begin{matrix} 
	-\frac{1}{x^2(x-2)^2} & \frac{x-1}{x^3(x-2)^2} \\
	\frac{2x-1}{x(x-2)^2} & -\frac{2x^2-4x+3}{x^2(x-2)^2}
	\end{matrix}},
\end{equation*}
then $L\psi=\psi F$ with 
\begin{equation*}
F(z)=\paren{\begin{matrix} 
	0 & 0 \\
	0 & z^2
	\end{matrix}}.
\end{equation*}
On the other hand, it  is easy to check that $\psi B=\theta \psi$ for 
\begin{equation*}
B=\partial_{z}^3.\paren{\begin{matrix} 
                   0 & 0 \\
	               1 & 0
	\end{matrix}}+ \partial_{z}^2.
\paren{\begin{matrix} 
	0 & 0 \\
	-\frac{2z+1}{z} & 0
	\end{matrix}}+
	\partial_{z}.\paren{\begin{matrix} 
		1 & 0 \\
		\frac{2(z-1)}{z^2} & 1
		\end{matrix}}
		+\paren{\begin{matrix} 
			-z^{-1} & 0 \\
			6z^{-3} & z^{-1}
			\end{matrix}}
\end{equation*}
and 
\begin{equation*}
\theta(x)=
\paren{\begin{matrix} 
	x & 0 \\
	x^2(x-2) & x
	\end{matrix}}.
\end{equation*}
In this opportunity we characterizes the algebra $\mathbb{A}$ of all polynomial $F$ such that there exist $L=L(x,\partial_{x})$ with $L\psi=\psi F$.

\begin{teor}[An algebra with two integral elements over one nilpotent and one idempotent]  \label{calogero}
	Let $\Gamma$ be the sub-algebra of $M_{2}(\complex)[z]$ of the form 
	\begin{equation*}
	\paren{\begin{matrix} 
		a & 0 \\
		b-a & b
		\end{matrix}}+
	\paren{\begin{matrix} 
		c & c \\
		a-b-c & -c
		\end{matrix}}z+
	\paren{\begin{matrix} 
		a-b-c & c+a-b\\
		d & e
		\end{matrix}}\frac{z^2}{2}+z^3p(z),
	\end{equation*}
	where $p\in M_{2}(\complex)[z]$ and all the variables $a,b,c,d,e$ are arbitrary. Then $\Gamma=\mathbb{A}$.
	
		Furthermore, we have the presentation $\mathbb{A}=\complex\cdot \langle \theta_{1}, \theta_{3},\theta_{4}, \theta_{5}\mid I=0\rangle$ with
	\begin{equation*}
	I=\langle\theta_{1}^2-\theta_{1},\theta_{4}^2,\theta_{4}\theta_{5},
	\theta_{4}\theta_{1}+\theta_{4}\theta_{3}-2\theta_{4}-\theta_{5}\theta_{4}-\theta_{5}^2,
	\theta_{3}^2-\theta_{3}+\theta_{5}-3\theta_{3}\theta_{4}\theta_{3}\theta_{5}-\theta_{1}\theta_{4}-\theta_{5}\theta_{1},
		\end{equation*}
	\begin{equation*}
	\theta_{3}\theta_{1}-\theta_{1}-\theta_{4}-\frac{1}{2}\theta_{4}\theta_{1}+\frac{1}{2}\theta_{4}\theta_{3}+\theta_{5}\theta_{1}-\frac{1}{2}\theta_{5}\theta_{4}
	+\frac{1}{2}\theta_{5}^2+\theta_{3}\theta_{4}-\theta_{1}\theta_{5}-\theta_{3}\theta_{5},
		\end{equation*}
	\begin{equation*}
	\theta_{1}\theta_{3}-\theta_{3}+\theta_{4}+\theta_{5}-\frac{3}{2}\theta_{4}\theta_{1}+\frac{3}{2}\theta_{4}\theta_{3}-2\theta_{5}\theta_{1}-\frac{3}{2}\theta_{5}\theta_{4}
	+\frac{3}{2}\theta_{5}^2+3\theta_{3}\theta_{4}+\theta_{3}\theta_{5},
		\end{equation*}
	\begin{equation*}
	\theta_{5}\theta_{3}-\theta_{4}\theta_{1}+\theta_{4}\theta_{3}-\theta_{5}\theta_{1}-\theta_{5}\theta_{4}+\theta_{5}^2,
	\theta_{5}\theta_{1}\theta_{5}-\theta_{5}^2\theta_{1}-\theta_{5}\theta_{4},
	\theta_{5}\theta_{4}\theta_{1}-\theta_{5}^3+\theta_{5}\theta_{1}\theta_{4}+\theta_{5}^2\theta_{1},
		\end{equation*}
	\begin{equation*}
	\theta_{4}\theta_{1}\theta_{5}+\theta_{4}\theta_{3}\theta_{5}-\theta_{3}^3,
	\theta_{5}\theta_{3}\theta_{4}+\theta_{5}\theta_{1}\theta_{4}
	\rangle.
	\end{equation*}
	
\end{teor}

This algebra is linked to the spin Calogero  systems whose relation with bispectrality can be found in \cite{BGK09}. See also~\cite{zbMATH06722531}. 

Theorems \ref{first1}, \ref{second} and \ref{calogero} give positive answers to the Conjectures 1, 2 and 3  of \cite{Grfrm-4} about three bispectral full rank 1 algebras. Moreover, these algebras are Noetherian and finitely generated because they are contained in the $N\times N$ matrix polynomial ring $M_{N}(\mathbb{K}[x])$. For more details see \cite{VZ1}.

The plan of this article is as follows: In Section \ref{sec1}, we consider noncommutative finitely generated algebras which are countably generated as left module over a subalgebra and prove Theorem \ref{presentation} about their presentations. In Section \ref{sec2}, we give positive answer to the first conjecture in \cite{Grfrm-4} by applying Theorem \ref{presentation} to obtain the presentation for an algebra with an integral element over a nilpotent one. In Section \ref{sec3}, we give positive answer to the second conjecture in \cite{Grfrm-4} by applying Theorem \ref{presentation} to obtain the presentation for an algebra with nilpotent and idempotent associated elements. Finally, in Section \ref{sec4}, we give positive answer to the third conjecture in \cite{Grfrm-4} by applying Theorem \ref{presentation} to obtain the presentation for an algebra with two integral elements over one nilpotent and one idempotent. 

\section{Presentations for Finitely Generated Algebras}\label{sec1}

In this section we face the presentation problem and obtain a method to tackle it. This method was motivated by a result used in the work presented in \cite{zurrian}. We begin by defining what a presentation is.
\begin{defin}
Let  $\mathbb{K}$ be a field, we denote by $\mathbb{K}\langle x_{\lambda}\mid \lambda\in \Lambda  \rangle$ the free algebra generated by the letters $x_{\lambda}$, $\lambda\in \Lambda$ i.e., $$\mathbb{K}\langle x_{\lambda}\mid \lambda\in \Lambda  \rangle=
\bigoplus_{F\subset \Lambda, F \text{finite}} \bigoplus_{\lambda\in F}\mathbb{K}\cdot x_{\lambda}.$$
\end{defin}
\begin{defin}
        Let $A$ be a $\mathbb{K}$-algebra. A presentation for an algebra $A$ is a triple $(\mathbb{K}\langle x_{\lambda}\mid \lambda\in \Lambda  \rangle,f,I)$ such that $I\subset A$ is an ideal and $f:\mathbb{K}\langle x_{\lambda}\mid \lambda\in \Lambda  \rangle/ I \rightarrow A$ is an isomorphism. Furthermore, we say that $A$ is finitely generated if there exists a presentation with $\Lambda$ finite and finitely presented if there exists a presentation with $\Lambda$ finite and the ideal $I$ is generated by finitely many elements.
\end{defin}

Now we are ready to state the main theorem.

\begin{teor}[Presentation of finitely generated algebras]\label{presentation}
	Let $A$ be a finitely generated $\mathbb{K}$-algebra by $\beta_{1}, \beta_{2},... , \beta_{n}$ such that:
	\begin{itemize}
			\item There exist an ideal  $I$  of $\mathbb{K}\cdot \langle \alpha_{1}, \alpha_{2}, ... , \alpha_{n} \rangle$ and an epimorphism of algebras
		\begin{equation*}
		f:\mathbb{K}\cdot \langle \alpha_{1}, \alpha_{2}, ... , \alpha_{n} \rangle /I \longrightarrow A,
		\end{equation*}
		\begin{equation*}
		f(\overline{\alpha_{j}})=\beta_{j}
		\end{equation*}

			\item There exists a subalgebra $\mathbb{K}\subset R\subset \mathbb{K}\cdot \langle \alpha_{1}, \alpha_{2}, ... , \alpha_{n} \rangle /I$ such that $\mathbb{K}\cdot \langle \alpha_{1}, \alpha_{2}, ... , \alpha_{n} \rangle /I$ is a free left $R$-module generated by $\llave{x_j}_{j=0}^{\infty}$, i.e.,  
		\begin{equation*}
		\mathbb{K}\cdot \langle \alpha_{1}, \alpha_{2}, ... , \alpha_{n} \rangle /I=\bigoplus_{j=0}^{\infty} Rx_j \mbox{ .}
		\end{equation*}  
		
		\item   $f\mid _{R}: R\longrightarrow A$ is a monomorphism.
		\item  The set $\llave{f(x_j)}_{j=0}^{\infty}$ is a basis for $A$ as a left $f(R)$-module.
		
		Then, $f$ is an isomorphism. 
		
	\end{itemize}
\end{teor}

\begin{demos}
It is enough to prove that $f$ is injective. Pick $x\in \ker(f)$ and write $x=\sum_{j=0}^{m}r_{j}x_j$, then $0=f(x)=\sum_{j=0}^{m}f(r_{j})f(x_j)$. However, since $\llave{f(x_j)}_{j=0}^{\infty}$ is a basis for $A$ as a left $f(R)$-module we have $f(r_{j})=0$ for $0\leq j \leq m$. Here we use that  $f\mid _{R}: R\longrightarrow A$ is an monomorphism to conclude $r_{j}=0$ for $0\leq j \leq m$ and $x=0$. \qed

\end{demos}
\begin{remark}
	The theorem guarantees a presentation of $A$ in terms of generators and relations through the isomorphism $f$, i.e., 
	\begin{equation*}
	A= \mathbb{K}\cdot \langle \beta_{1}, \beta_{2}, ... , \beta_{n}\mid P(\beta_{1}, \beta_{2}, ... , \beta_{n})=0, \forall P\in I \rangle .
	\end{equation*}
\end{remark}
This theorem is a method to find out presentations for finitely generated algebras. However, we need to choose generators for the  algebra and look for relations among them. Furthermore, we must seek for an intermediate $\mathbb{K}$-algebra $\mathbb{K}\subset S\subset A$ and a linearly independent set $\llave{y_j}_{j=0}^{\infty}$ of $A$ such that $A=\bigoplus_{j=0}^{\infty} Sy_j$.

In the following sections we shall apply this method to obtain presentations for some noncommutative bispectral algebras.

\section{An algebra with an integral element over a nilpotent one}\label{sec2}

In this section, we give the proof of Theorem \ref{first1}. The idea of the proof is to consider a basis for the vector space 
$\mathbb{A}\cap \bigoplus_{j=0}^{3} M_{2}(\complex[x])_{j}$ of polynomials in $\mathbb{A}$ of degree least or equal to $3$ and observe that this basis generates the algebra $\mathbb{A}$. After that, we look for remarkable elements in the basis that generate the others and obtain some set of relations. Finally, we verify the hypothesis of the Theorem \ref{presentation} to obtain the proof of the assertion.

	
\begin{demos}[Theorem~\ref{first1}]


	The first part of the proof is given by the Theorem 3 in \cite{VZ1}. We will give a proof of the existence of the presentation.
Note that $\mathbb{A}$ is generated by $\beta_{0}=e_{12}$, $\beta_{1}=Ix+e_{21}x^2$, $\beta_{2}=e_{12}x+e_{11}x^2$, $\beta_{3}=e_{12}x+e_{22}x^2$, 
	$\beta_{4}=e_{12}x^2$, $\beta_{5}=e_{12}x-e_{21}x^3$, $\beta_{6}=e_{11}x^3$, $\beta_{7}=e_{12}x^3$, $\beta_{8}=e_{22}x^3$.
	
	Moreover, we can eliminate the variables $\beta_j$ for $2\leq j\leq 8$. In fact, 
	$\beta_2=\beta_{0}\beta_{1}$, $\beta_3=\beta_{1}\beta_{0}$, $\beta_4=\beta_{0}\beta_{1}\beta_{0}$,
	$\beta_5=\frac{\beta_{0}\beta_{1}+\beta_{1}\beta_{0}-\beta_{1}^2}{2}$, $\beta_6=\frac{\beta_{0}\beta_{1}\beta_{0}-\beta_{0}\beta_{1}^2}{2}$,
	$\beta_7=\frac{\beta_{0}\beta_{1}^2\beta_{0}}{2}$,
	 $\beta_8=\frac{\beta_{0}\beta_{1}\beta_{0}-\beta_{1}^2\beta_{0}}{2}$.

Furthermore, we are going to check the presentation using Theorem \ref{presentation}. We begin with some general results:

\begin{prop}\label{propspin}
	Let $A$ be  a $\mathbb{K}$-algebra. Suppose that $\beta_0 \in A$ is a nilpotent element of degree $2$, then 
	\begin{equation*}
	\llave{\beta_{1}^j \mid j\geq 0} \cup \llave{\beta_{1}^j \beta_{0} \mid j\geq 0} 
	\cup\llave{\beta_{1}^j \beta_{0} \beta_{1} \mid j\geq 0} \cup \llave{\beta_{1}^j \beta_{0} \beta_{1} \beta_{0}\mid j\geq 0} 
	\end{equation*}
	is a linearly independent set over $\mathbb{K}$ if and only if 
	\begin{equation*}
\llave{\beta_{1}^j \beta_{0} \mid j\geq 0} \cup \llave{\beta_{1}^j \beta_{0} \beta_{1} \beta_{0}\mid j\geq 0} 
	\end{equation*}
	is a linearly independent set over $\mathbb{K}$.
\end{prop}
\begin{demos}
	Clearly the condition is sufficient. We consider the expression:
	\begin{equation}\label{ind 2x2}
	\sum_{j=0}^{n}a_{j}\beta_{1}^j+	\sum_{j=0}^{n}b_{j}\beta_{1}^j \beta_{0}+ 	\sum_{j=0}^{n}c_{j}\beta_{1}^j \beta_{0}\beta_{1}
		+\sum_{j=0}^{n}d_{j}\beta_{1}^j \beta_{0}\beta_{1}\beta_{0}=0
	\end{equation}
	for $a_j, b_j, c_j, d_j \in \mathbb{K}$, $n\in \natu$.
	
Multiply by $\beta_{0}$ on the right and using that $\beta_{0}^2=0$ we obtain:
	\begin{equation*}
\sum_{j=0}^{n}a_{j}\beta_{1}^j \beta_{0}+ \sum_{j=0}^{n}c_{j}\beta_{1}^j \beta_{0}\beta_{1}\beta_{0}=0.
\end{equation*}
If we assume that $\llave{\beta_{1}^j \beta_{0} \mid j\geq 0} \cup \llave{\beta_{1}^j \beta_{0} \beta_{1} \beta_{0}\mid j\geq 0} $ is linearly independent we have $a_j = c_j =0$ and 
\eqref{ind 2x2} reduces to:
	\begin{equation*}
\sum_{j=0}^{n}b_{j}\beta_{1}^j \beta_{0}+ \sum_{j=0}^{n}d_{j}\beta_{1}^j \beta_{0}\beta_{1}\beta_{0}=0.
\end{equation*}
Again, using this assumption we have  $b_j = d_j =0$. With this fact we obtain the necessity. \qed
\end{demos}

\begin{prop}\label{ind2x2 prop}
	Taking the elements $\beta_{0}$ and $\beta_{1}$ in $\mathbb{A}$ we obtain that
		\begin{equation*}
	\llave{\beta_{1}^j \beta_{0} \mid j\geq 0} \cup \llave{\beta_{1}^j \beta_{0} \beta_{1} \beta_{0}\mid j\geq 0} 
	\end{equation*}
	is a linearly independent set.
\end{prop}
\begin{demos}
	Note that $\beta_{1}^j \beta_{0}=e_{12}x^j+je_{22}x^{j+1}$ and $\beta_{1}^j \beta_{0} \beta_{1} \beta_{0}=e_{12}x^{j+2}+je_{22}x^{j+3}$. Consider the expression:
		\begin{equation*}
	\sum_{j=0}^{n}a_{j}\beta_{1}^j \beta_{0}+ \sum_{j=0}^{n}b_{j}\beta_{1}^j \beta_{0}\beta_{1}\beta_{0}=0.
	\end{equation*}
	Replacing the previous relations we obtain:
		\begin{equation*}
	\sum_{j=0}^{n}a_{j}(e_{12}x^j+je_{22}x^{j+1})+ \sum_{j=0}^{n}b_{j}\beta_{1}^j (e_{12}x^{j+2}+je_{22}x^{j+3})=0.
	\end{equation*}
	
Using the entries of the matrix we obtain:
\begin{equation*}
\sum_{j=0}^{n}a_j x^j+ \sum_{j=0}^{n}b_j x^{j+2}=0  \hspace{0.1cm} \text{and } \sum_{j=0}^{n}j a_j x^{j+1}+ \sum_{j=0}^{n}jb_j x^{j+3}=0.  
\end{equation*}
Equivalently,
\begin{equation*}
\sum_{j=0}^{n}a_j x^j+ \sum_{j=2}^{n+2}b_{j-2} x^{j}=0  \hspace{0.1cm} \text{and } \sum_{j=0}^{n}j a_j x^{j}+ \sum_{j=2}^{n+2}(j-2)b_{j-2} x^{j}=0.  
\end{equation*}
Hence, 
\begin{equation*}
a_0+a_1 x+\sum_{j=2}^{n}(a_j + b_{j-2}) x^j+ b_{n-1}x^{n+1}+b_{n}x^{n+2}=0  \hspace{0.1cm} \text{and } 
\end{equation*}
\begin{equation*}
a_1 x+ \sum_{j=2}^{n}(j a_j+ (j-2)b_{j-2}) x^{j}+(n-1)b_{n-1}x^{n+1}+nb_{n}x^{n+2} =0.  
\end{equation*}
Therefore, 
\begin{equation*}
a_0 =a_1 = b_{n-1}=b_{n}=0,
\paren{\begin{matrix} 
	1 & 1\\
	j & j-2
	\end{matrix}}
\paren{\begin{matrix} 
	a_j \\
	b_{j-2}
	\end{matrix}}=\paren{\begin{matrix} 
	0\\
	0
	\end{matrix}},  2\leq j \leq n.
\end{equation*}
Since $\det \paren{\begin{matrix} 
	1 & 1\\
	j & j-2
	\end{matrix}} = -2\neq 0$ we have $a_j=b_{j-2}=0$,\hspace{0.1cm} $2\leq j \leq n$ and 
	\begin{equation*}
\llave{\beta_{1}^j \beta_{0} \mid j\geq 0} \cup \llave{\beta_{1}^j \beta_{0} \beta_{1} \beta_{0}\mid j\geq 0} 
\end{equation*}
is  linearly independent.\qed
\end{demos}
\begin{lemma}\label{lemmaspin}
	Consider the algebra $\mathbb{K} \cdot \langle \alpha_{0}, \alpha_{1} \rangle/I$ with 
	\begin{equation*}
	I=<\alpha_{0}^2,\alpha_{1}^3+\alpha_{0}\alpha_{1}\alpha_{0}-3\alpha_{1}\alpha_{0}\alpha_{1}+\alpha_{0}\alpha_{1}^2+\alpha_{1}^2\alpha_{0}>
	\end{equation*}
	then 	$\llave{1, \alpha_{0}, \alpha_{0} \alpha_{1}, \alpha_{0} \alpha_{1} \alpha_{0}}   $ 
	 is a system of generators for $\mathbb{K}\cdot \langle \alpha_{0}, \alpha_{1} \rangle/I$ as a free left $R$-module,
	with $R=\mathbb{K}\cdot \langle \alpha_{1} \rangle /I$.
\end{lemma}
\begin{demos}
	Define $M=R  \oplus R \cdot \alpha_{0} \oplus R \cdot \alpha_{0} \alpha_{1}  \oplus R \cdot \alpha_{0} \alpha_{1} \alpha_{0}$. We have to see that $\mathbb{K}\cdot \langle \alpha_{0}, \alpha_{1} \rangle/I=M$. It is enough to show that $M$ is invariant under left and right multiplication by $\alpha_{0}$ and $\alpha_{1}$.
	\begin{itemize}
		\item $\alpha_1 M \subset M$. 
		
		Since $\alpha_{1} \in R$.
		
		\item  $M \alpha_0 \subset M$.  In fact,
$M \alpha_{0} \subset 
	R \cdot \alpha_{0} \oplus   R \cdot \alpha_{0}\alpha_{1}\alpha_{0} \subset M$.
		
		\item  $M \alpha_{1} \subset M$. 
		
	Since $\alpha_{0}\alpha_{1}^2=-\alpha_{1}^3-\alpha_{1}^2 \alpha_{0} +3\alpha_{1}\alpha_{0}\alpha_{1}-\alpha_{0}\alpha_{1}\alpha_{0}$ we have
	$$\alpha_{0}\alpha_{1}^2 \alpha_{0} =-\alpha_{1}^3 \alpha_{0} +3(\alpha_{1}\alpha_{0})^2$$ and $$0= -\alpha_{0}\alpha_{1}^3\alpha_{0}-\alpha_{0}\alpha_{1}^2\alpha_{0}
	+3(\alpha_{0}\alpha_{1})^2.$$
	
	Furthermore, 
	\begin{equation*}
	\alpha_{0}\alpha_{1}^3=-\alpha_{1}^4-\alpha_{1}^2(\alpha_{0}\alpha_{1})+3\alpha_{1}(\alpha_{0}\alpha_{1}^2)-(\alpha_{0}\alpha_{1})^2.
	\end{equation*}
		Hence,
		\begin{equation*}
		3(\alpha_{0}\alpha_{1})^2=\alpha_{0}\alpha_{1}^3+
		\alpha_{0}\alpha_{1}^2 \alpha_{0}=	\alpha_{0}\alpha_{1}^3-\alpha_{1}^3 \alpha_{0} +3(\alpha_{1}\alpha_{0})^2
		\end{equation*}
			\begin{equation*}
	=-\alpha_{1}^4-\alpha_{1}^2(\alpha_{0}\alpha_{1})+3\alpha_{1}(\alpha_{0}\alpha_{1}^2)-(\alpha_{0}\alpha_{1})^2-\alpha_{1}^3 \alpha_{0} +3(\alpha_{1}\alpha_{0})^2.
		\end{equation*}
		
Equivalently,
	\begin{equation*}
	4(\alpha_{0}\alpha_{1})^2=-\alpha_{1}^4-\alpha_{1}^3 \alpha_{0}-\alpha_{1}^2(\alpha_{0}\alpha_{1})+3\alpha_{1}(\alpha_{0}\alpha_{1}^2)+3(\alpha_{1}\alpha_{0})^2.
\end{equation*}
		
However,
	\begin{equation*}
\alpha_{1}\alpha_{0}\alpha_{1}^2=-\alpha_{1}^4-\alpha_{1}^3 \alpha_{0}+3\alpha_{1}^2(\alpha_{0}\alpha_{1})-(\alpha_{1}\alpha_{0})^2.
\end{equation*}

Thus,
	\begin{equation*}
4(\alpha_{0}\alpha_{1})^2=-\alpha_{1}^4-\alpha_{1}^3 \alpha_{0}-\alpha_{1}^2(\alpha_{0}\alpha_{1})+3\alpha_{1}^4-3\alpha_{1}^3\alpha_{0}+9\alpha_{1}^2(\alpha_{0}\alpha_{1})-3(\alpha_{1}\alpha_{0})^2
+3(\alpha_{1}\alpha_{0})^2=
\end{equation*}
	\begin{equation*}
-4\alpha_{1}^4-4\alpha_{1}^3\alpha_{0}+8\alpha_{1}^2 (\alpha_{0}\alpha_{1}).
\end{equation*}
Therefore, 
	\begin{equation*}
(\alpha_{0}\alpha_{1})^2=
-\alpha_{1}^4-\alpha_{1}^3\alpha_{0}+2\alpha_{1}^2 (\alpha_{0}\alpha_{1}).
\end{equation*}
This implies that $(\alpha_{0}\alpha_{1})^2 \in M$, $\alpha_{0}\alpha_{1}^2 \in M$. Since $M$ is a left $R$-module we have
$M\alpha_{1}\subset R\alpha_{1} \oplus R\alpha_{0}\alpha_{1} \oplus  R\alpha_{0}\alpha_{1}^2 \oplus  R(\alpha_{0}\alpha_{1})^2 \subset M$.
		\item $ \alpha_{0} M \subset M$. 
		
		We claim that $\alpha_{0}\alpha_{1}^n \in M$ for every $n\in \natu$. 
		For $n=0$ is clear. Assume this for some $n\in \natu$ and note that $\alpha_{0}\alpha_{1}^{n+1}
		=(\alpha_{0}\alpha_{1}^n)\alpha_{1}\in M\alpha_{1}\subset M$. The claim follows by induction.
		
		In particular, $\alpha_{0}R\subset M$. Thus, $\alpha_{0}M\subset \alpha_{0} R\oplus \alpha_{0} R\alpha_{0} \oplus  \alpha_{0} R\alpha_{0}\alpha_{1}
		 \oplus \alpha_{0}  R \alpha_{0}\alpha_{1}\alpha_{0} \subset R  \oplus R \cdot \alpha_{0} \oplus R \cdot \alpha_{0} \alpha_{1}  \oplus R \cdot \alpha_{0} \alpha_{1} \alpha_{0}
		 \subset  M$.

\qed

Finally, we conclude with the proof of the nice presentation.  
Define 
\begin{equation*}
f:\complex\cdot \langle \alpha_{0}, \alpha_{1}\rangle /I \longrightarrow \mathbb{A},
\end{equation*}
\begin{equation*}
f(\overline{\alpha_{j}})=\beta_{j}
\end{equation*}
the previous lemma guarantees the existence of a subalgebra $R=\complex\cdot \langle \alpha_{1} \rangle /I$ and  a system of generators
$\llave{1,\alpha_{0}, \alpha_{0}\alpha_{1},\alpha_{0}\alpha_{1}\alpha_{0}}$   for $\complex \cdot \langle \alpha_{0}, \alpha_{1} \rangle/I$ as a free left $R$-module. Furthermore,
$f\mid _{R}: R\longrightarrow A$ is a monomorphism.

The Proposition \ref{ind2x2 prop} implies that $\llave{1,\beta_{0}, \beta_{0}\beta_{1},\beta_{0}\beta_{1}\beta_{0}}$ is a linearly independent set over $\complex$. Thus, we are under the hypothesis of Theorem \ref{presentation} and $f$ is an isomorphism.  \qed
	\end{itemize}
\end{demos}
Putting together Lemma~\ref{lemmaspin}, Propositions~\ref{propspin} and \ref{ind2x2 prop}, we conclude the proof of Theorem~\ref{first1}.  \qed
\end{demos}

\section{An algebra with nilpotent and idempotent associated elements}\label{sec3}

In this section, we give the proof of Theorem \ref{first1}. The idea of the proof is to consider a basis for the vector space 
$\mathbb{A}\cap \bigoplus_{j=0}^{5} M_{3}(\complex[x])_{j}$ of polynomials in $\mathbb{A}$ of degree least or equal to $5$ and observe that this basis generates the algebra $\mathbb{A}$. After that, we look for remarkable elements in the basis that generate the others and obtain some set of relations. Finally, we verify the hypothesis of the Theorem \ref{presentation} to obtain the proof of the assertion.


 
 \begin{demos}
    The proof is a straightforward check of the relations given in the Proposition 1 in \cite{VZ1}. We will give a proof of the presentation. 
    
Note that $\mathbb{A}$ is generated by 	$\beta_{0}=e_{13},	\beta_{1}=e_{12}-e_{33}x+e_{21}x^2+e_{31}x^3, \beta_{2}=e_{12}+e_{23},$\\
	$\beta_{3}=e_{22}+(e_{21}+e_{32})x+e_{31}x^2,	\beta_{4}=e_{22}x+S_{3}x^2+e_{31}x^3,	\beta_{5}=Ix-e_{31}x^3,	\beta_{6}=e_{13}x-e_{11}x^3,$ \\	$\beta_{7}=e_{13}x-e_{22}x^3,
	\beta_{8}=e_{13}x-e_{33}x^3,	\beta_{9}=e_{13}x^2,	\beta_{10}=e_{13}x-e_{31}x^5,	\beta_{11}=e_{23}x-e_{33}x^2,$ \\	$\beta_{12}=S_{3}x+e_{31}x^4,	\beta_{13}=Ix^2-2e_{31}x^4,	\beta_{14}=e_{12}x^2+e_{31}x^5,	\beta_{15}=e_{22}x^2-e_{31}x^4,	\beta_{16}=e_{23}x^2-e_{31}x^5,$ 	\\ $ \beta_{17}=e_{12}x^3,
	\beta_{18}=e_{13}x^3,	\beta_{19}=e_{21}x^3+e_{31}x^4,	\beta_{20}=e_{23}x^3,	\beta_{21}=e_{32}x^3+e_{31}x^4,	\beta_{22}=e_{11}x^4,
	\beta_{23}=e_{12}x^4,$ \\ $	\beta_{24}=e_{13}x^4,	\beta_{25}=e_{21}x^4+e_{31}x^5,	\beta_{26}=e_{22}x^4,	\beta_{27}=e_{23}x^4,	\beta_{28}=e_{32}x^4+e_{31}x^5,
	\beta_{29}=e_{33}x^4,	\beta_{30}=e_{11}x^5,$ \\	$ \beta_{31}=e_{12}x^5,	\beta_{32}=e_{13}x^5,	\beta_{33}=e_{21}x^5,	\beta_{34}=e_{22}x^5,
	\beta_{35}=e_{23}x^5,	\beta_{36}=e_{32}x^5,	\beta_{37}=e_{33}x^5$. \\
	However, we can eliminate the variables $\beta_{j}$ for $j\neq 2,3$. In fact, \\
		$\beta_{0}=\beta_{2}^2, \hspace{0.1cm}  \beta_{1}=1/2\beta_{3}\beta_{2}\beta_{3}-\beta_{3}\beta_{2}+\beta_{2}, 
	\hspace{0.1cm} \beta_{4}=	1/2\beta_{3}\beta_{2}\beta_{3}, \hspace{0.1cm}
	\beta_{5}=-1/2\beta_{3}\beta_{2}\beta_{3}+\beta_{2}\beta_{3}+\beta_{3}\beta_{2}-\beta_{2},$\\
	$\beta_{6}=	-1/2\beta_{2}^2\beta_{3}\beta_{2}\beta_{3}+\beta_{2}^2\beta_{3}\beta_{2},
	\hspace{0.1cm}  \beta_{7}=	-1/2\beta_{2}(\beta_{3}\beta_{2})^2+\beta_{2}^2\beta_{3}\beta_{2}+\beta_{2}\beta_{3}\beta_{2}^2,
	\hspace{0.1cm}\beta_{8}=	-1/2(\beta_{3}\beta_{2})^2\beta_{2}+\beta_{2}\beta_{3}\beta_{2}^2,
\\ \beta_{9}=\beta_{2}^2\beta_{3}\beta_{2}^2, 
\beta_{10}=	-1/2(\beta_{3}\beta_{2})^2\beta_{2}\beta_{3}-1/2\beta_{2}(\beta_{2}\beta_{3})^2+\beta_{2}\beta_{3}\beta_{2}^2\beta_{3}-1/2\beta_{2}(\beta_{3}\beta_{2})^2+\beta_{3}\beta_{2}^2\beta_{3}\beta_{2}\\-1/2(\beta_{3}\beta_{2})^2\beta_{2}+\beta_{2}^2\beta_{3}\beta_{2}
+\beta_{2}\beta_{3}\beta_{2}^2, \hspace{0.1cm} 
\beta_{11}=	\beta_{3}\beta_{2}^2,
	\hspace{0.1cm} \beta_{12}=	-1/2(\beta_{2}\beta_{3})^2+\beta_{3}\beta_{2}^2\beta_{3}-1/2(\beta_{3}\beta_{2})^2+\beta_{2}^2\beta_{3}+\beta_{2}\beta_{3}\beta_{2}+\beta_{3}\beta_{2}^2-\beta_{2}^2,
	\hspace{0.1cm}  \beta_{13}=	(\beta_{2}\beta_{3})^2-2\beta_{3}\beta_{2}^2\beta_{3}+(\beta_{3}\beta_{2})^2-\beta_{2}^2\beta_{3}-\beta_{2}\beta_{3}\beta_{2}-\beta_{3}\beta_{2}^2+\beta_{2}^2, 
	\hspace{0.1cm}  \beta_{14}=	1/2\beta_{3}\beta_{2}\beta_{3}\beta_{2}\beta_{2}\beta_{3}+1/2\beta_{2}\beta_{2}\beta_{3}\beta_{2}\beta_{3}-\beta_{2}\beta_{3}\beta_{2}\beta_{2}\beta_{3}+1/2\beta_{2}\beta_{3}\beta_{2}\beta_{3}\beta_{2}-\beta_{3}\beta_{2}\beta_{2}\beta_{3}\beta_{2}+1/2\beta_{3}\beta_{2}\beta_{3}\beta_{2}\beta_{2}-\beta_{2}\beta_{3}\beta_{2}\beta_{2},
	\hspace{0.1cm}  \beta_{15}=	1/2(\beta_{2}\beta_{3})^2-\beta_{3}\beta_{2}^2\beta_{3}+1/2(\beta_{3}\beta_{2})^2-\beta_{2}^2\beta_{3}-\beta_{3}\beta_{2}^2,
	\hspace{0.1cm}  \beta_{16}= 1/2\beta_{3}\beta_{2}\beta_{3}\beta_{2}^2\beta_{3}+1/2\beta_{2}(\beta_{2}\beta_{3})^2-\beta_{2}\beta_{3}\beta_{2}^2\beta_{3}+1/2\beta_{2}(\beta_{3}\beta_{2})^2-\beta_{3}\beta_{2}^2\beta_{3}\beta_{2}+1/2(\beta_{3}\beta_{2})^2\beta_{2}-\beta_{2}^2\beta_{3}\beta_{2},
	\hspace{0.1cm}  \beta_{22}=	\beta_{2}^2\beta_{3}\beta_{2}^2\beta_{3}-1/2\beta_{2}^2(\beta_{3}\beta_{2})^2+\beta_{2}^2\beta_{3}\beta_{2}^2,
	\hspace{0.1cm}  \beta_{17}=	1/2\beta_{2}^2(\beta_{3}\beta_{2})^2-\beta_{2}^2\beta_{3}\beta_{2}^2,
	\hspace{0.1cm} \beta_{18}=	1/2\beta_{2}(\beta_{2}\beta_{3})^2\beta_{2}^2,
	\hspace{0.1cm}  \beta_{19}=	\beta_{3}\beta_{2}^2\beta_{3}-1/2(\beta_{3}\beta_{2})^2+\beta_{3}\beta_{2}^2,
	\hspace{0.1cm}  \beta_{29}=	-1/2(\beta_{2}\beta_{3})^2\beta_{2}^2+\beta_{3}\beta_{2}^2\beta_{3}\beta_{2}^2+\beta_{2}^2\beta_{3}\beta_{2}^2,
	\hspace{0.1cm}  \beta_{20}=	1/2(\beta_{2}\beta_{3})^2\beta_{2}^2-\beta_{2}^2\beta_{3}\beta_{2}^2,
	\hspace{0.1cm}  \beta_{21}=	-1/2(\beta_{2}\beta_{3})^2+\beta_{3}\beta_{2}^2\beta_{3}+\beta_{2}^2\beta_{3},
	\hspace{0.1cm}  \beta_{23}=\beta_{2}^2\beta_{3}\beta_{2}^2\beta_{3}\beta_{2}-1/2\beta_{2}^2\beta_{3}\beta_{2}\beta_{3}\beta_{2}^2,
	\hspace{0.1cm}  \beta_{24}=-\beta_{2}^2\beta_{3}\beta_{2}^2\beta_{3}\beta_{2}^2,
	\hspace{0.1cm}  \beta_{25}=	1/2\beta_{3}\beta_{2}\beta_{3}\beta_{2}^2\beta_{3}-\beta_{3}\beta_{2}^2\beta_{3}\beta_{2}+1/2\beta_{3}\beta_{2}\beta_{3}\beta_{2}^2,
	\hspace{0.1cm}  \beta_{26}=	-1/2\beta_{2}^2\beta_{3}\beta_{2}\beta_{3}\beta_{2}+\beta_{2}\beta_{3}\beta_{2}^2\beta_{3}\beta_{2}-1/2\beta_{2}\beta_{3}\beta_{2}\beta_{3}\beta_{2}^2+\beta_{2}^2\beta_{3}\beta_{2}^2,
	\hspace{0.1cm}  \beta_{37}=	1/2\beta_{3}\beta_{2}\beta_{3}\beta_{2}^2\beta_{3}\beta_{2}^2+1/2\beta_{2}^2\beta_{3}\beta_{2}\beta_{3}\beta_{2}^2-\beta_{2}\beta_{3}\beta_{2}^2\beta_{3}\beta_{2}^2,
	\hspace{0.1cm}  \beta_{27}=-1/2\beta_{2}^2\beta_{3}\beta_{2}\beta_{3}\beta_{2}^2+\beta_{2}\beta_{3}\beta_{2}^2\beta_{3}\beta_{2}^2,
	\hspace{0.1cm}  \beta_{28}=	1/2\beta_{3}\beta_{2}\beta_{3}\beta_{2}^2\beta_{3}+1/2\beta_{2}^2\beta_{3}\beta_{2}\beta_{3}-\beta_{2}\beta_{3}\beta_{2}^2\beta_{3},
	\hspace{0.1cm}  \beta_{30}	=1/2\beta_{2}^2\beta_{3}\beta_{2}\beta_{3}\beta_{2}^2\beta_{3}-\beta_{2}^2\beta_{3}\beta_{2}^2\beta_{3}\beta_{2}+1/2\beta_{2}^2\beta_{3}\beta_{2}\beta_{3}\beta_{2}^2,
	\hspace{0.1cm}  \beta_{36}=	1/2\beta_{3}\beta_{2}\beta_{3}\beta_{2}^2\beta_{3}\beta_{2}+1/2\beta_{2}^2(\beta_{3}\beta_{2})^2-\beta_{2}\beta_{3}\beta_{2}^2\beta_{3}\beta_{2}+1/2(\beta_{2}\beta_{3})^2\beta_{2}^2-\beta_{3}\beta_{2}^2\beta_{3}\beta_{2}^2-\beta_{2}^2\beta_{3}\beta_{2}^2,
	\hspace{0.1cm}  \beta_{31}=	1/2\beta_{2}^2\beta_{3}\beta_{2}\beta_{3}\beta_{2}^2\beta_{3}\beta_{2}-\beta_{2}^2\beta_{3}\beta_{2}^2\beta_{3}\beta_{2}^2,
	\hspace{0.1cm}  \beta_{32}=	1/2\beta_{2}^2\beta_{3}\beta_{2}\beta_{3}\beta_{2}^2\beta_{3}\beta_{2}^2,
	\hspace{0.1cm}  \beta_{33}=	1/2(\beta_{2}\beta_{3})^2\beta_{2}^2\beta_{3}-\beta_{2}^2\beta_{3}\beta_{2}^2\beta_{3}+1/2\beta_{2}^2(\beta_{3}\beta_{2})^2-\beta_{2}\beta_{3}\beta_{2}^2\beta_{3}\beta_{2}
	+1/2(\beta_{2}\beta_{3})^2\beta_{2}^2-\beta_{2}^2\beta_{3}\beta_{2}^2,
	\hspace{0.1cm}  \beta_{34}= 1/2(\beta_{2}\beta_{3})^2\beta_{2}^2\beta_{3}\beta_{2}-\beta_{2}^2\beta_{3}\beta_{2}^2\beta_{3}\beta_{2}+1/2\beta_{2}^2(\beta_{3}\beta_{2})^2\beta_{2}-\beta_{2}\beta_{3}\beta_{2}^2\beta_{3}\beta_{2}^2, \\ \beta_{35}=	1/2\beta_{2}\beta_{3}\beta_{2}\beta_{3}\beta_{2}^2\beta_{3}\beta_{2}^2-\beta_{2}^2\beta_{3}\beta_{2}^2\beta_{3}\beta_{2}^2$.
Furthermore, we are going to check the presentation using  Theorem \ref{presentation}. We begin with some general results:
\begin{lemma}\label{indk}
	Let $A$ be  a $\mathbb{K}$-algebra. Suppose that $\beta_2 \in A$ is a nilpotent element of degree $D\geq 3$. Suppose that 
	\begin{equation*}
	\llave{\beta_{2}^{D-1}(\beta_{3}\beta_{2})^j \beta_{2}^{D-2} \mid  j\geq 0}
	\end{equation*}
	is a linearly independent set over $\mathbb{K}$. Then, 
	$\llave{\beta_{2}^{D-1}(\beta_{3}\beta_{2})^j \beta_{2}^{k} \mid  j\geq 0, 1\leq k \leq D-2}$ is linearly independent over $\mathbb{K}$.
\end{lemma}
\begin{demos}	
	Consider the expression 
	\begin{equation}\label{ind3}
	\sum_{j=1}^{n}\sum_{k=1}^{D-2}c_{jk}\beta_{2}^{D-2} (\beta_3 \beta_2)^j \beta_{2}^k=0.
	\end{equation}
	Multiplying by $\beta_{2}^{D-3}$ on the right:
	\begin{equation}
	\sum_{j=1}^{n}\sum_{k=1}^{D-2}c_{j1}\beta_{2}^{D-2} (\beta_3 \beta_2)^j \beta_{2}^{D-2}=0.
	\end{equation}
	However, $\llave{\beta_{2}^{D-1}(\beta_{3}\beta_{2})^j \beta_{2}^{D-2} \mid  j\geq 0}$ is linearly independent over $\mathbb{K}$. Thus $c_{j1}=0$ for $0\leq j \leq n$.
	
	Thus \eqref{ind3} reduces to 
	\begin{equation}
	\sum_{j=1}^{n}\sum_{k=2}^{D-2}c_{jk}\beta_{2}^{D-2} (\beta_3 \beta_2)^j \beta_{2}^k=0.
	\end{equation}
	Assume that 
	\begin{equation}
	\sum_{j=1}^{n}\sum_{k=k_0}^{D-2}c_{jk}\beta_{2}^{D-2} (\beta_3 \beta_2)^j \beta_{2}^k=0.
	\end{equation}
	Multiplying by $\beta_{2}^{D-2-k_0}$ on the right:
	\begin{equation}
	\sum_{j=1}^{n}c_{jk_0}\beta_{2}^{D-2} (\beta_3 \beta_2)^j \beta_{2}^{D-2}=0.
	\end{equation}
	However, $\llave{\beta_{2}^{D-1}(\beta_{3}\beta_{2})^j \beta_{2}^{D-2} \mid  j\geq 0}$ is linearly independent over $k$. Thus $c_{jk_0}=0$ for $1\leq j \leq n$. Thus
	\begin{equation}
	\sum_{j=1}^{n}\sum_{k=k_{0}+1}^{D-2}c_{jk}\beta_{2}^{D-2} (\beta_3 \beta_2)^j \beta_{2}^k=0.
	\end{equation}
	Since the case $k_0 =1\Rightarrow k_0 =2$ was seen we have that $c_{jk}=0$ for $1\leq j \leq n, 1\leq k \leq D-2$. \qed
\end{demos}
\begin{prop}\label{ind3x3}
Let $A$ be  a $\mathbb{K}$-algebra. Suppose that $\beta_2 \in A$ is a nilpotent element of degree $D\geq 3$, then 
\begin{equation*}
\llave{\beta_{2}^i (\beta_{3}\beta_{2})^j \beta_3 \mid 0\leq i\leq D-1, j\geq 0}\cup \llave{\beta_{2}^i (\beta_{3}\beta_{2})^j \beta_{2}^{k} \mid 0\leq i\leq D-1, j\geq 1, 1\leq k \leq D-2}
\end{equation*}
\begin{equation*}
\cup \llave{\beta_{2}^i (\beta_{3}\beta_{2})^j \mid 0\leq i\leq D-1, j\geq 0}
\end{equation*}
is a linearly independent set over $\mathbb{K}$ if and only if 
\begin{equation*}
\llave{\beta_{2}^{D-1}(\beta_{3}\beta_{2})^j \beta_{2}^{D-2} \mid  j\geq 0}
\end{equation*}
is a linearly independent set over $\mathbb{K}$.
\end{prop}

\begin{demos}
	The sufficiency of the  statement is clear. To show the necessity we consider the expresion 
	\begin{equation}\label{ind1}
	\sum_{j=0}^{n}\sum_{i=0}^{D-1}a_{ij}\beta_{2}^{i} (\beta_3 \beta_2)^j + 	\sum_{j=0}^{n}\sum_{i=0}^{D-1}b_{ij}\beta_{2}^{i} (\beta_3 \beta_2)^j \beta_3
	+	\sum_{j=1}^{n}\sum_{i=0}^{D-1}\sum_{k=1}^{D-2}c_{ijk}\beta_{2}^{i} (\beta_3 \beta_2)^j \beta_{2}^k=0
	\end{equation}
	$a_{ij},b_{ij},c_{ijk}\in \mathbb{K}, n\geq 0$.
	
	We have to see that 	$a_{ij}=b_{ij}=c_{ijk}=0$.

	We are going to see that 
		\begin{equation}\label{indl}
	\sum_{j=0}^{n}\sum_{i=l}^{D-1}a_{ij}\beta_{2}^{i} (\beta_3 \beta_2)^j + 	\sum_{j=0}^{n}\sum_{i=l}^{D-1}b_{ij}\beta_{2}^{i} (\beta_3 \beta_2)^j \beta_3
	+	\sum_{j=1}^{n}\sum_{i=l}^{D-1}\sum_{k=1}^{D-2}c_{ijk}\beta_{2}^{i} (\beta_3 \beta_2)^j \beta_{2}^k=0
	\end{equation}
	for some $0\leq l\leq D-1$ implies that $a_{lj}=b_{lj}=c_{ljk}=0$.

For $l=0$ we have the equation \eqref{ind1}. Multiplying by $\beta_{2}^{D-1}$ on the left and on the right:
\begin{equation}
	\sum_{j=0}^{n}\sum_{i=0}^{D-1}b_{ij}\beta_{2}^{D-1} (\beta_3 \beta_2)^j \beta_3 \beta_{2}^{D-1}=
	\sum_{j=0}^{n}\sum_{i=0}^{D-1}b_{ij}\beta_{2}^{D-1} (\beta_3 \beta_2)^{j+1} \beta_{2}^{D-2}=0.
\end{equation}
However, $\llave{\beta_{2}^{D-1}(\beta_{3}\beta_{2})^j \beta_{2}^{D-2} \mid  j\geq 0}$ is linearly independent over $\mathbb{K}$. Thus, $b_{0j}=0$ for $0\leq j \leq n$.

This reduces \eqref{ind1} to 
	\begin{equation}
\sum_{j=0}^{n}\sum_{i=0}^{D-1}a_{ij}\beta_{2}^{i} (\beta_3 \beta_2)^j  +	\sum_{j=1}^{n}\sum_{i=0}^{D-1}\sum_{k=1}^{D-2}c_{ijk}\beta_{2}^{i} (\beta_3 \beta_2)^j \beta_{2}^k=0.
\end{equation}
Multiplying by $\beta_{2}^{D-1}$ on the left:
	\begin{equation}\label{ind2}
\sum_{j=0}^{n}a_{0j}\beta_{2}^{D-1} (\beta_3 \beta_2)^j  +	\sum_{j=1}^{n}\sum_{k=1}^{D-2}c_{0jk}\beta_{2}^{D-2} (\beta_3 \beta_2)^j \beta_{2}^k=0.
\end{equation}
Multiplying by  $\beta_{2}^{D-2}$ on the right:
	\begin{equation}
\sum_{j=0}^{n}a_{0j}\beta_{2}^{D-1} (\beta_3 \beta_2)^j \beta_{2}^{D-2}  =0.
\end{equation}
Thus, $a_{0j}=0$ for $0\leq j \leq n$. Since $\llave{\beta_{2}^{D-1}(\beta_{3}\beta_{2})^j \beta_{2}^{D-2} \mid  j\geq 0}$ is linearly independent over $\mathbb{K}$. 

This reduces \eqref{ind2} to 
	\begin{equation}
\sum_{j=0}^{n}\sum_{k=1}^{D-2}c_{0jk}\beta_{2}^{D-2} (\beta_3 \beta_2)^j \beta_{2}^k=0.
\end{equation}

However, by Lemma \ref{indk}, $\llave{\beta_{2}^{D-1}(\beta_{3}\beta_{2})^j \beta_{2}^{k} \mid  j\geq 0, 1\leq k \leq D-2}$ is linearly independent over $\mathbb{K}$. Thus $c_{0jk}=0$ for $1\leq j \leq n, 1\leq k \leq D-2$.

Assume \eqref{indl} for $l$ and multiply this by $\beta_{2}^{D-l-1}$ on the left:
	\begin{equation}\label{indl1}
\sum_{j=0}^{n}a_{lj}\beta_{2}^{D-1} (\beta_3 \beta_2)^j + 	\sum_{j=0}^{n}b_{lj}\beta_{2}^{D-1} (\beta_3 \beta_2)^j \beta_3
+	\sum_{j=1}^{n}\sum_{k=1}^{D-2}c_{ljk}\beta_{2}^{D-1} (\beta_3 \beta_2)^j \beta_{2}^k=0.
\end{equation}
Multiplying by $\beta_{2}^{D-1}$ on the right:
	\begin{equation}
	\sum_{j=1}^{n}b_{lj}\beta_{2}^{D-1} (\beta_3 \beta_2)^j \beta_3 \beta_{2}^{D-1}
=\sum_{j=1}^{n}\sum_{k=1}^{D-2}b_{lj}\beta_{2}^{D-1} (\beta_3 \beta_2)^{j+1} \beta_{2}^{D-2}=0.
\end{equation}
However, $\llave{\beta_{2}^{D-1}(\beta_{3}\beta_{2})^j \beta_{2}^{D-2} \mid  j\geq 0}$ is linearly independent over $\mathbb{K}$. Thus, $b_{lj}=0$ for $0\leq j \leq n$.

Therefore, \eqref{indl1} reduces to:
	\begin{equation}
\sum_{j=0}^{n}a_{lj}\beta_{2}^{D-1} (\beta_3 \beta_2)^j 
+	\sum_{j=1}^{n}\sum_{k=1}^{D-2}c_{ljk}\beta_{2}^{D-1} (\beta_3 \beta_2)^j \beta_{2}^k=0.
\end{equation}
Multiplying by $\beta_{2}^{D-2}$ on the right:
	\begin{equation}
\sum_{j=0}^{n}a_{lj}\beta_{2}^{D-1} (\beta_3 \beta_2)^j \beta_{2}^{D-2}=0.
\end{equation}
However, $\llave{\beta_{2}^{D-1}(\beta_{3}\beta_{2})^j \beta_{2}^{D-2} \mid  j\geq 0}$ is linearly independent over $\mathbb{K}$. Thus, $a_{lj}=0$ for $0\leq j \leq n$.

Therefore,
\begin{equation} 
\label{indl1bis}
	\sum_{j=1}^{n}\sum_{k=1}^{D-2}c_{ljk}\beta_{2}^{D-1} (\beta_3 \beta_2)^j \beta_{2}^k=0.
\end{equation}
However, by Lemma \ref{indk}, \hspace{0.1 cm} $\llave{\beta_{2}^{D-1}(\beta_{3}\beta_{2})^j \beta_{2}^{k} \mid  j\geq 0, 1\leq k \leq D-2}$ is linearly independent over $\mathbb{K}$. Thus, $c_{ljk}=0$ for $1\leq j \leq n, 1\leq k \leq D-2$.

Thus, we obtain \eqref{indl1bis} for $l+1$. Then  \eqref{indl1bis}  is valid for $0\leq l \leq D-1$, i.e., $a_{ij}=b_{ij}=c_{ijk}=0$. \qed
\end{demos}

\begin{lemma}\label{lemma3x3}
Consider the algebra $\mathbb{K}\cdot \langle \alpha_{2}, \alpha_{3} \rangle/I$ with 
\begin{equation*}
I=<\alpha_{2}^3,\alpha_{3}^2-\alpha_{3}, (\alpha_{3}\alpha_{2})^2 \alpha_{3}-4\alpha_{3}\alpha_{2}^2 \alpha_{3}>
\end{equation*}
then $\llave{(\alpha_{3}\alpha_{2})^n \mid n\geq 0}\cup \llave{(\alpha_{3}\alpha_{2})^n  \alpha_{3}\mid n\geq 0}\cup
\llave{(\alpha_{3}\alpha_{2})^n  \alpha_{2}\mid n\geq 0}$ is a system of generators for $\mathbb{K}\cdot \langle \alpha_{2}, \alpha_{3} \rangle/I$ as a free left $R$-module,
with $R=\mathbb{K}\cdot \langle \alpha_{2} \rangle /I$.
\end{lemma}
\begin{demos}
Define $M=\bigoplus_{n=0}^{\infty} R \cdot (\alpha_{3}\alpha_{2})^n \oplus 
\bigoplus_{n=0}^{\infty} R \cdot (\alpha_{3}\alpha_{2})^n \alpha_{3} \oplus
\bigoplus_{n=1}^{\infty} R \cdot (\alpha_{3}\alpha_{2})^n \alpha_{2}$. We have to see that $\mathbb{K}\cdot \langle \alpha_{2}, \alpha_{3} \rangle/I=M$. It is enough to show that $M$ is invariant under left and right multiplication by $\alpha_{2}$ and $\alpha_{3}$.
\begin{itemize}
	\item $\alpha_2 M \subset M$. 
	
	Since $\alpha_{2} \in R$.
	
	\item  $M \alpha_2 \subset M$. 
	
	Since $R (\alpha_{3}\alpha_{2})^n \alpha_{2}\subset M$, $\corch{(\alpha_{3}\alpha_{2})^n \alpha_{3}}\alpha_{2}=
	(\alpha_{3}\alpha_{2})^{n+1}\in M$, for $n\geq 0$ and $\corch{(\alpha_{3}\alpha_{2})^n \alpha_{2}}\alpha_{2}=0 \in M$ for $n\geq 1$.
	Then $M \alpha_{2} \subset 
	\bigoplus_{n=0}^{\infty} R \cdot (\alpha_{3}\alpha_{2})^{n+1} \oplus 
	\bigoplus_{n=0}^{\infty} R \cdot (\alpha_{3}\alpha_{2})^n \alpha_{2}\subset M$.
	
	\item  $M \alpha_{3} \subset M$. 
	
	Note that  $\corch{(\alpha_{3}\alpha_{2})^n \alpha_{2}}\alpha_{3}=(\alpha_{3}\alpha_{2})^n \alpha_{2}\alpha_{3}=
	(\alpha_{3}\alpha_{2})^{n-1} \alpha_{3}\alpha_{2}^2 \alpha_{3}=
	\frac{1}{4}	(\alpha_{3}\alpha_{2})^{n-1} (\alpha_{3}\alpha_{2})^2 \alpha_{3}
	=\frac{1}{4}	(\alpha_{3}\alpha_{2})^{n+1} \alpha_{3}$ for every $n\geq 1$, then 
	$M\alpha_{3} \subset 
	\bigoplus_{n=0}^{\infty} R \cdot (\alpha_{3}\alpha_{2})^n \alpha_{3} \oplus 
	\bigoplus_{n=0}^{\infty} R \cdot (\alpha_{3}\alpha_{2})^n \alpha_{3} \oplus
	\bigoplus_{n=1}^{\infty} R \cdot (\alpha_{3}\alpha_{2})^n \alpha_{2} \alpha_{3}
	\subset \bigoplus_{n=0}^{\infty} R \cdot (\alpha_{3}\alpha_{2})^n \alpha_{3} \oplus 
	\bigoplus_{n=1}^{\infty} R \cdot (\alpha_{3}\alpha_{2})^{n+1} \alpha_{3} \subset M$.
	
	\item  $M \alpha_{3} \subset M$. 
	
	Note that  $\corch{(\alpha_{3}\alpha_{2})^n \alpha_{2}}\alpha_{3}=(\alpha_{3}\alpha_{2})^n \alpha_{2}\alpha_{3}=
	(\alpha_{3}\alpha_{2})^{n-1} \alpha_{3}\alpha_{2}^2 \alpha_{3}=
	\frac{1}{4}	(\alpha_{3}\alpha_{2})^{n-1} (\alpha_{3}\alpha_{2})^2 \alpha_{3}
	=\frac{1}{4}	(\alpha_{3}\alpha_{2})^{n+1} \alpha_{3}$ for every $n\geq 1$, then 
	$M\alpha_{3} \subset 
	\bigoplus_{n=0}^{\infty} R \cdot (\alpha_{3}\alpha_{2})^n \alpha_{3} \oplus 
	\bigoplus_{n=0}^{\infty} R \cdot (\alpha_{3}\alpha_{2})^n \alpha_{3} \oplus
	\bigoplus_{n=1}^{\infty} R \cdot (\alpha_{3}\alpha_{2})^n \alpha_{2} \alpha_{3}
	\subset \bigoplus_{n=0}^{\infty} R \cdot (\alpha_{3}\alpha_{2})^n \alpha_{3} \oplus 
	\bigoplus_{n=1}^{\infty} R \cdot (\alpha_{3}\alpha_{2})^{n+1} \alpha_{3} \subset M$.
	
	\item $ \alpha_{3} M \subset M$. 
	
	Note that $$\alpha_{3}\alpha_{2}^2 (\alpha_{3}\alpha_{2})^n =
	(\alpha_{3}\alpha_{2}^2 \alpha_{3}) \alpha_{2} (\alpha_{3}\alpha_{2})^{n-1}
	=\frac{1}{4}\corch{ (\alpha_{3}\alpha_{2})^2 \alpha_{3}}\alpha_{2}  (\alpha_{3}\alpha_{2})^{n-1}
	=\frac{1}{4} (\alpha_{3}\alpha_{2})^2 \alpha_{3}\alpha_{2}  (\alpha_{3}\alpha_{2})^{n-1}$$
	$=\frac{1}{4} (\alpha_{3}\alpha_{2})^{n+2}\in M$ for $n\geq 1$ and $\alpha_{3}\alpha_{2}^2 =
	(\alpha_{3}\alpha_{2})\alpha_{2}\in M$. Then  $\alpha_{3}\alpha_{2}^2 (\alpha_{3}\alpha_{2})^n \in M$ for every $n\geq 0$. 
	
	On the other hand $\alpha_{3}\alpha_{2} (\alpha_{3}\alpha_{2})^n =(\alpha_{3}\alpha_{2})^{n+1}\in M$, $\alpha_{3}  (\alpha_{3}\alpha_{2})^n=
	 (\alpha_{3}\alpha_{2})^n \in M$ for all  $n\geq 0$.
	 
	 Furthermore, $$\alpha_{3} (\alpha_{3}\alpha_{2})^n \alpha_3 =  (\alpha_{3}\alpha_{2})^n \alpha_{3} \in M, 
	 \alpha_{3} (\alpha_{3}\alpha_{2})^n \alpha_2 = (\alpha_{3}\alpha_{2})^n \alpha_{2} \in M,$$ for all $n\geq 0$, and
	 $$(\alpha_{3} \alpha_{2}) (\alpha_{3}\alpha_{2})^n \alpha_3 =  (\alpha_{3}\alpha_{2})^{n+1}\alpha_{3} \in M, 	
	 (\alpha_{3} \alpha_{2}) (\alpha_{3}\alpha_{2})^n \alpha_2 =(\alpha_{3} \alpha_{2})^{n+1} \alpha_{2} \in M$$ for all $n\geq 0$.
	 
	 On the other hand $(\alpha_{3}\alpha_{2}^2) (\alpha_{3}\alpha_{2})^n \alpha_{3}= \frac{1}{4} (\alpha_{3}\alpha_{2})^{n+2} \alpha_{3}\in M$,
	  $(\alpha_{3}\alpha_{2}^2) (\alpha_{3}\alpha_{2})^n \alpha_{2}= \frac{1}{4} (\alpha_{3}\alpha_{2})^{n+2} \alpha_{2}\in M$ for all $n\geq 0$. In particular $\alpha_{3} M\subset M$. \qed
\end{itemize}
\end{demos}

Finally, we conclude with the proof of the nice presentation.  
Define 
	\begin{equation*}
f:\complex\cdot \langle \alpha_{2}, \alpha_{3}\rangle /I \longrightarrow \mathbb{A},
\end{equation*}
\begin{equation*}
f(\overline{\alpha_{j}})=\beta_{j}
\end{equation*}
the previous lemma guarantees the existence of a subalgebra $R=\complex\cdot \langle \alpha_{2} \rangle /I$ and  a system of generators
$\llave{(\alpha_{3}\alpha_{2})^n \mid n\geq 0}\cup \llave{(\alpha_{3}\alpha_{2})^n  \alpha_{3}\mid n\geq 0}\cup
\llave{(\alpha_{3}\alpha_{2})^n  \alpha_{2}\mid n\geq 0}$   for $\complex \cdot \langle \alpha_{2}, \alpha_{3} \rangle/I$ as a free left $R$-module. Furthermore
 $f\mid _{R}: R\longrightarrow A$ is a monomorphism.
 
Since $\beta_{2}^2 (\beta_3 \beta_2)^n \beta_2 =2^{n-1} e_{13}x^{n+1}$ for $n\geq 1$ applying the Proposition \ref{ind3x3} with $D=3$ we obtain 
\begin{equation*}
\llave{\beta_{2}^i (\beta_{3}\beta_{2})^j \beta_3 \mid 0\leq i\leq D-1, j\geq 0}\cup \llave{\beta_{2}^i (\beta_{3}\beta_{2})^j \beta_{2} \mid 0\leq i\leq D-1, j\geq 1}
\end{equation*}
\begin{equation*}
\cup \llave{\beta_{2}^i (\beta_{3}\beta_{2})^j \mid 0\leq i\leq D-1, j\geq 0}
\end{equation*}
is a linearly independent set over $\complex$. 

	Putting together  Proposition~\ref{ind3x3}, Lemmas~\ref{indk} and \ref{lemma3x3} we conclude the proof of Theorem~\ref{second}. \qed
 \end{demos}

\section{An algebra with two integral elements over one nilpotent and one idempotent}\label{sec4}

In this section we give the proof of Theorem~\ref{calogero}.  The idea of the proof is to consider a basis for the vector space 
$\mathbb{A}\cap \bigoplus_{j=0}^{2} M_{2}(\complex[x])_{j}$ of polynomials in $\mathbb{A}$ of degree least or equal to $2$ and observe that this basis generates the algebra $\mathbb{A}$. After that, we look for remarkable elements in the basis that generate the others and obtain some set of relations. Finally, we verify the hypothesis of the Theorem \ref{presentation} to obtain the proof of the assertion.



\begin{demos}

The first part of the proof is given by the Theorem 3 in \cite{VZ1}. We are going to check the presentation using  Theorem \ref{presentation}.

\begin{lemma}\label{calogen}
	Consider the $\mathbb{K}$-algebra $\mathbb{K}\cdot \langle \theta_{1}, \theta_{3}, \theta_{4}, \theta_{5} \rangle/I$ with $\mathbb{K}$ a central field of characteristic $0$ and 
	\begin{equation*}
I=\langle\theta_{1}^2-\theta_{1},\theta_{4}^2,\theta_{4}\theta_{5},
\theta_{4}\theta_{1}+\theta_{4}\theta_{3}-2\theta_{4}-\theta_{5}\theta_{4}-\theta_{5}^2,
\theta_{3}^2-\theta_{3}+\theta_{5}-3\theta_{3}\theta_{4}-\theta_{3}\theta_{5}-\theta_{1}\theta_{4}-\theta_{5}\theta_{1},
\end{equation*}
\begin{equation*}
\theta_{3}\theta_{1}-\theta_{1}-\theta_{4}-\frac{1}{2}\theta_{4}\theta_{1}+\frac{1}{2}\theta_{4}\theta_{3}+\theta_{5}\theta_{1}-\frac{1}{2}\theta_{5}\theta_{4}
+\frac{1}{2}\theta_{5}^2+\theta_{3}\theta_{4}-\theta_{1}\theta_{5}-\theta_{3}\theta_{5},
\end{equation*}
\begin{equation*}
\theta_{1}\theta_{3}-\theta_{3}+\theta_{4}+\theta_{5}-\frac{3}{2}\theta_{4}\theta_{1}+\frac{3}{2}\theta_{4}\theta_{3}-2\theta_{5}\theta_{1}-\frac{3}{2}\theta_{5}\theta_{4}
+\frac{3}{2}\theta_{5}^2+3\theta_{3}\theta_{4}+\theta_{3}\theta_{5},
\end{equation*}
\begin{equation*}
\theta_{5}\theta_{3}-\theta_{4}\theta_{1}+\theta_{4}\theta_{3}-\theta_{5}\theta_{1}-\theta_{5}\theta_{4}+\theta_{5}^2,
\theta_{5}\theta_{1}\theta_{5}-\theta_{5}^2\theta_{1}-\theta_{5}\theta_{4},
\theta_{5}\theta_{4}\theta_{1}-\theta_{5}^3+\theta_{5}\theta_{1}\theta_{4}+\theta_{5}^2\theta_{1},
\end{equation*}
\begin{equation*}
\theta_{4}\theta_{1}\theta_{5}+\theta_{5}^2\theta_{1}+\theta_{5}\theta_{4}-\theta_{5}^3,
\theta_{5}\theta_{3}\theta_{4}+\theta_{5}\theta_{1}\theta_{4},
 \theta_{3}\theta_{5} \theta_{1}-\theta_{1}\theta_{5}-\theta_{3}\theta_{5}+\theta_{3}\theta_{4}+\theta_{5}^2,
 \end{equation*}
 \begin{equation*}
 \theta_{4}\theta_{1}\theta_{5}+\theta_{4}\theta_{3}\theta_{5}-\theta_{5}^3,
 \theta_{4}\theta_{3}\theta_{5}-\theta_{5}\theta_{1}\theta_{5},
 \theta_{5}\theta_{3}\theta_{4}+\theta_{5}\theta_{1}\theta_{4},
  \end{equation*}
 \begin{equation*}
\theta_{1}\theta_{5}\theta_{1}+\theta_{1}\theta_{5}+\theta_{3}\theta_{5}-\theta_{5}^2+\theta_{1}\theta_{4}-\theta_{5}^2\theta_{1}-\theta_{5}\theta_{4}
-\theta_{3}\theta_{5}^2,
  \end{equation*}
\begin{equation*}
 \theta_{1}\theta_{5}\theta_{3}-2\theta_{3}\theta_{4}-\theta_{1}\theta_{5}-\theta_{3}\theta_{5}+\theta_{5}^2+3\theta_{5}\theta_{4}-
 \theta_{1}\theta_{4}+\theta_{5}^2\theta_{1}+\theta_{3}\theta_{5}^2+2\theta_{5}\theta_{1}\theta_{4}+2\theta_{3}\theta_{5}\theta_{4},
   \end{equation*}
 \begin{equation*}
 \theta_{3}\theta_{4}\theta_{1}+\theta_{1}\theta_{4}+\theta_{1}\theta_{5}+\theta_{3}\theta_{5}-\theta_{5}\theta_{4}-\theta_{3}\theta_{5}^2-\theta_{5}^2 \rangle.
\end{equation*}

	Then, $\llave{\theta_{4}\theta_{1},\theta_{3},\theta_{1}}\cup \llave{\theta_{5}^n\mid n\geq 0}\cup \llave{\theta_{5}^n \theta_{4} \mid n\geq 0}\cup
	\llave{\theta_{5}^n \theta_{1}\theta_{4}  \mid n\geq 0}\cup \llave{\theta_{5}^n \theta_{1} \mid n\geq 1}\cup
\llave{\theta_{3}\theta_{5}^n \mid n\geq 1}\cup \llave{\theta_{1}\theta_{5}^n \mid  n\geq 1}\cup \llave{\theta_{3}\theta_{5}^n\theta_{4} \mid n\geq 0}
\cup \llave{\theta_{1}\theta_{5}^n\theta_{4} \mid n\geq 1}$ is a system of generators for $\mathbb{K}\cdot \langle \theta_{1}, \theta_{3}, \theta_{4}, \theta_{5} \rangle/I$  as a free $\mathbb{K}$-vector space.
\end{lemma}
\begin{demos}
	Define $M=\mathbb{K} \cdot \theta_{1} \oplus \mathbb{K} \cdot \theta_{3} \oplus \mathbb{K}\cdot\theta_{4}\theta_{1} \oplus 
	\bigoplus_{n=0}^{\infty} \mathbb{K} \cdot \theta_{5}^n\oplus
	\bigoplus_{n=0}^{\infty} \mathbb{K} \cdot\theta_{5}^n \theta_{4} \oplus 
	\bigoplus_{n=0}^{\infty} \mathbb{K}\cdot \theta_{5}^n \theta_{1}\theta_{4} \oplus 
	\bigoplus_{n=1}^{\infty} \mathbb{K}\cdot \theta_{5}^n \theta_{1} \oplus 
	\bigoplus_{n=1}^{\infty} \mathbb{K} \cdot \theta_{3}\theta_{5}^n \oplus 
	\bigoplus_{n=1}^{\infty} \mathbb{K} \cdot \theta_{1}\theta_{5}^n \oplus 
	\bigoplus_{n=0}^{\infty} \mathbb{K} \cdot \theta_{3}\theta_{5}^n\theta_{4}\oplus 
	\bigoplus_{n=1}^{\infty} \mathbb{K} \cdot \theta_{1}\theta_{5}^n\theta_{4} $. We have to see that $\mathbb{K}\cdot \langle \theta_{1}, \theta_{3}, \theta_{4}, \theta_{5} \rangle/I=M$. It is enough to show that $M$ is invariant under left and right multiplication by $\theta_{1},\theta_{3},\theta_{4}$ and $\theta_{5}$.
	
	\begin{itemize}
		\item $M\theta_{5} \subset M$. 
		
		Note that $\theta_{3}\theta_{4}\theta_{5}=0\in M$, $\theta_{4}\theta_{1}\theta_{5}=-\theta_{5}^2\theta_{1}-\theta_{5}\theta_{4}+\theta_{5}^3\in M$. On the other hand 
		$\theta_{1}\theta_{4}\theta_{5}=0\in M$, $\theta_{3}\theta_{5}\in M$, $\theta_{1}\theta_{5}\in M$, $\theta_{5}^n \in M$ for every $n\geq 1$, $\theta_{5}^n \theta_{4}\theta_{5}=0
		\in M$, for every $n\geq 0$, $\theta_{5}^n \theta_{1}\theta_{4}\theta_{5}=0\in M$, for every $n\geq 1$.
		
		Furthermore, $\theta_{5}^n \theta_{1}\theta_{5}= \theta_{5}^{n+1}\theta_{1}+\theta_{5}^n \theta_{4}\in M$ for every $n\geq 1$, $(\theta_{3}\theta_{5}^n)\theta_{5}
		=\theta_{3}\theta_{5}^{n+1}\in M$ for every $n\geq 1$, $(\theta_{1}\theta_{5}^n)\theta_{5}=\theta_{1}\theta_{5}^{n+1}\in M$ for every $n\geq 1$, 
		$(\theta_{3}\theta_{5}^n \theta_{4})\theta_{5}=0\in M$ for every $n\geq 1$, 	$(\theta_{1}\theta_{5}^n \theta_{4})\theta_{5}=0\in M$ for every $n\geq 1$. In particular 
		$M\theta_{5} \subset M$. 
		
		\item  $M\theta_{4} \subset M$. 
		
		Note that $\theta_{4}\in M$, $(\theta_{3}\theta_{4})\theta_{4}=0\in M$, $\theta_{4}\theta_{1}=\theta_{5}\theta_{3}+\theta_{4}\theta_{3}-\theta_{5}\theta_{1}-\theta_{5}
		\theta_{4}+\theta_{5}^2$ and $\theta_{4}\theta_{3}=2\theta_{4}+\theta_{5}\theta_{4}+\theta_{5}^2-\theta_{4}\theta_{1}$ imply $\theta_{4}\theta_{1}=
		\frac{1}{2}\theta_{5}\theta_{3}+\theta_{4}+\theta_{5}^2-\frac{1}{2}\theta_{5}\theta_{1}$, hence $\theta_{4}\theta_{1}\theta_{4}=
		\theta_{5}^2 \theta_{4}-\theta_{5}\theta_{1}\theta_{4} \in M$. 
		
		On the other hand, $(\theta_{1} \theta_{4})\theta_{4}=\theta_{1}\theta_{4}^2=0\in M$,  $\theta_{3}\theta_{4}\in M$, $\theta_{1}\theta_{4}\in M$, $\theta_{5}^n \theta_{4}
		\in M$ for every $n\geq 0$, $(\theta_{5}^n \theta_{4})\theta_{4}=0\in M$ for every $n\geq 0$, $(\theta_{5}^n \theta_{1}\theta_{4})\theta_{4}=0\in M$ for every 
		$n\geq 1$, $(\theta_{5}^n \theta_{1})\theta_{4}=\theta_{5}^n \theta_{1}\theta_{4}\in M$ for every $n\geq 1$, $(\theta_{3}\theta_{5}^n)\theta_{4}=
		\theta_{3}\theta_{5}^n \theta_{4}\in M$ for every $n\geq 1$, $(\theta_{1}\theta_{5}^n)\theta_{4}=
		\theta_{1}\theta_{5}^n \theta_{4}\in M$ for every $n\geq 1$, $(\theta_{3}\theta_{5}^n \theta_{4})\theta_{4}=0\in M$ for every $n\geq 1$, 
		$(\theta_{1}\theta_{5}^n \theta_{4})\theta_{4}=0\in M$ for every $n\geq 1$. In particular  $M\theta_{4} \subset M$. 
		
		\item  $\theta_{1} M \subset M$.
		
		Note that $\theta_{1}\in M$. Since $\theta_{1}\theta_{3}=\theta_{3}-\theta_{4}-\theta_{5}+\frac{3}{2}\theta_{4}\theta_{1}-\frac{3}{2}\theta_{4}\theta_{3}+2\theta_{5}\theta_{1}
		+\frac{3}{2}\theta_{5}\theta_{4}-\frac{3}{2}\theta_{5}^2-3\theta_{3}\theta_{4}-\theta_{3}\theta_{5}$ multiplying by $\theta_{4}$ on the right we obtain $\theta_{1}\theta_{3}
		\theta_{4}=\theta_{3}\theta_{4}-\theta_{5}\theta_{4}-\theta_{5}\theta_{1}\theta_{4}-\theta_{3}\theta_{5}\theta_{4}\in M$.
		
		On the other hand, the equation $\theta_{4}\theta_{1}=
		\frac{1}{2}\theta_{5}\theta_{3}+\theta_{4}+\theta_{5}^2-\frac{1}{2}\theta_{5}\theta_{1}$ implies 
		$\theta_{1}\theta_{4}\theta_{1}=
		\frac{1}{2}\theta_{1}\theta_{5}\theta_{3}+\theta_{1}\theta_{4}+\theta_{1}\theta_{5}^2-\frac{1}{2}\theta_{1}\theta_{5}\theta_{1}$. Putting this equation  together with the equations 
		 \begin{equation*}
\theta_{1}\theta_{5}\theta_{1}+\theta_{1}\theta_{5}+\theta_{3}\theta_{5}-\theta_{5}^2+\theta_{1}\theta_{4}-\theta_{5}^2\theta_{1}-\theta_{5}\theta_{4}
-\theta_{3}\theta_{5}^2=0,
  \end{equation*}
\begin{equation*}
 \theta_{1}\theta_{5}\theta_{3}-2\theta_{3}\theta_{4}-\theta_{1}\theta_{5}-\theta_{3}\theta_{5}+\theta_{5}^2+3\theta_{5}\theta_{4}-
 \theta_{1}\theta_{4}+\theta_{5}^2\theta_{1}+\theta_{3}\theta_{5}^2+2\theta_{5}\theta_{1}\theta_{4}+2\theta_{3}\theta_{5}\theta_{4}=0,
   \end{equation*}
		we obtain:
		$$\theta_{1}\theta_{4}\theta_{1}-\theta_{3}\theta_{4}-2\theta_{1}\theta_{4}-\theta_{1}\theta_{5}-\theta_{3}\theta_{5}+2\theta_{5}\theta_{4}+\theta_{5}\theta_{1}\theta_{4}
		+\theta_{3}\theta_{5}\theta_{4}+\theta_{5}^2\theta_{1}+\theta_{3}\theta_{5}^2-\theta_{1}\theta_{5}^2+\theta_{5}^2=0.$$
		In particular,  $\theta_{1}\theta_{4}\theta_{1}\in M$.
		
		Moreover, 
		$\theta_{1}\theta_{5}\theta_{1}+\theta_{1}\theta_{4}+\theta_{1}\theta_{5}+\theta_{3}\theta_{5}-\theta_{5}^2\theta_{1}-\theta_{5}\theta_{4}-\theta_{3}\theta_{5}^2
		-\theta_{5}^2=0$ implies $\theta_{1}\theta_{5}\theta_{1}\in M$.
		
		However, multiplying  $\theta_{4}\theta_{1}+\theta_{4}\theta_{3}-2\theta_{4}-\theta_{5}\theta_{4}-\theta_{5}^2=0$ by $\theta_{1}$ on the left we have
		 $\theta_{1}\theta_{4}\theta_{1}+\theta_{1}\theta_{4}\theta_{3}-2\theta_{1}\theta_{4}-\theta_{1}\theta_{5}\theta_{4}-\theta_{1}\theta_{5}^2=0$. Hence,
		 $\theta_{1}\theta_{4}\theta_{3}=\theta_{5}^2-\theta_{3}\theta_{4}-\theta_{1}\theta_{5}-\theta_{3}\theta_{5}+2\theta_{5}\theta_{4}+\theta_{5}\theta_{1}\theta_{4}+\theta_{3}\theta_{5}\theta_{4}
		 +\theta_{5}^2\theta_{1}+\theta_{3}\theta_{5}^2+\theta_{1}\theta_{5}\theta_{4}\in M$.
		
		Moreover, $\theta_{1}(\theta_{1}\theta_{4})=\theta_{1}\theta_{4}\in M$ and $\theta_{1}\theta_{3}=\theta_{3}-\theta_{4}-\theta_{5}+\frac{3}{2}\theta_{4}\theta_{1}-\frac{3}{2}
		\theta_{4}\theta_{3}+2\theta_{5}\theta_{1}+\frac{3}{2}\theta_{5}\theta_{4}-\frac{3}{2}\theta_{5}^2-3\theta_{3}\theta_{4}-\theta_{3}\theta_{5}\in M$. On the other hand, 
		$\theta_{1}^2=\theta_{1}\in M$, $\theta_{5}^n\in M$ for every  $n\geq 1$, $\theta_{1}\theta_{5}^n \theta_{4}\in M$ for every $n\geq 0$. Note that $\theta_{1}\theta_{5}^n \theta_{1}
		=-\theta_{1}\theta_{5}^n-\theta_{3}\theta_{5}^n-\theta_{1}\theta_{5}^{n-1}\theta_{4}+\theta_{5}^{n+1}\theta_{1}+\theta_{5}^n \theta_{4}+\theta_{3}\theta_{5}^{n+1}+\theta_{5}^{n+1}\in M$ for every $n\geq 2$ and $\theta_{1}\theta_{5}\theta_{1}\in M$ imply $\theta_{1}\theta_{5}^n \theta_{1}\in M$ for every 
		$n\geq 1$. 
		
		Since $\theta_{1}\theta_{3}\in M$ we have $\theta_{1}\theta_{3}\theta_{5}^n\in M\theta_{5}^n\subset M$ for every $n\geq 1$. Furthermore, $\theta_{1}(\theta_{1}\theta_{5}^n)=
		\theta_{1}\theta_{5}^n\in M$ for every $n\geq 1$ and $\theta_{1}(\theta_{3}\theta_{5}^n\theta_{4})=(\theta_{1}\theta_{3}\theta_{5}^n)\theta_{4}\in M\theta_{4}\subset M$ for every 
			$n\geq 1$.  However, $\theta_{1}\theta_{3}\theta_{4}\in M$ then $\theta_{1}(\theta_{3}\theta_{5}^n\theta_{4})\in M$ for every $n\geq 0$. Note that $\theta_{1}(\theta_{1}\theta_{5}^n\theta_{4})=\theta_{1}\theta_{5}^n\theta_{4}\in M$ for every $n\geq 1$. Thus $\theta_{1}M\subset M$.
			
			\item $\theta_{4}M\subset M$.
			
			Note that $\theta_{4}(\theta_{4}\theta_{1})=0\in M$. Furthermore, $\theta_{4}\theta_{3}=2\theta_{4}+\theta_{5}\theta_{4}+\theta_{5}^2-\theta_{4}\theta_{1}\in M$ and $\theta_{4}\theta_{1}\in M$, $\theta_{4}\theta_{5}^n=0\in M$ for every $n\geq 1$ and $\theta_{4}\in M$. Moreover, $\theta_{4}\theta_{5}^n\theta_{4}=0\in M$ for every $n\geq 0$, 
			$\theta_{4}(\theta_{5}^n\theta_{1}\theta_{4})=0\in M$ for every $n\geq 1$. Since $\theta_{4}(\theta_{1}\theta_{4})=(\theta_{4}\theta_{1})\theta_{4}\in M\theta_{4}\subset M$ we have $\theta_{4}\theta_{5}^n \theta_{1}\theta_{4}\in M$ for every $n\geq 0$. 
			
			On the other hand, $\theta_{4}(\theta_{5}\theta_{1})=0\in M$ for every $n \geq 1$. Since $\theta_{4}\theta_{3}\in M$ we have $\theta_{4}(\theta_{3}\theta_{5}^n)\in  M\theta_{5}^n\subset M$ for every $n\geq 1$.  Since $\theta_{4}\theta_{1}$ we have $\theta_{4}(\theta_{1}\theta_{5}^n)\in M\theta_{5}^n\subset M$ for every $n\geq 1$. Using that $\theta_{4}\theta_{3}\theta_{5}^n\in M$ we have $\theta_{4}(\theta_{3}\theta_{5}^n\theta_{4})\in M\theta_{4}\subset M$. Since $\theta_{4}\theta_{1}\theta_{5}^n\in M$ we have 
			$\theta_{4}(\theta_{1}\theta_{5}^n\theta_{4})=(\theta_{4}\theta_{1}\theta_{5}^n)\theta_{4}\in M\theta_{4}\subset M$. Thus, $\theta_{4}M\subset M$.
			
			\item $M\theta_{1}\subset M$.
			
			Note that $(\theta_{4}\theta_{1})\theta_{1}=\theta_{4}\theta_{1}\in M$. Since $\theta_{3}\theta_{1}=\theta_{1}+\theta_{4}+\frac{1}{2}\theta_{4}\theta_{1}-\frac{1}{2}\theta_{4}\theta_{3}-\theta_{5}\theta_{1}+\frac{1}{2}\theta_{5}\theta_{4}-\frac{1}{2}\theta_{5}^2-\theta_{3}\theta_{4}+\theta_{1}\theta_{5}+\theta_{3}\theta_{5}$ and $\theta_{4}\theta_{3}=2\theta_{4}+\theta_{5}\theta_{4}+\theta_{5}^2-\theta_{4}\theta_{1}$ we have 
			$\theta_{3}\theta_{1}=\theta_{1}+\theta_{4}\theta_{1}-\theta_{5}^2-\theta_{5}\theta_{1}+\theta_{1}\theta_{5}+\theta_{3}\theta_{5}-\theta_{3}\theta_{4}\in M$.
			
			On the other hand, $\theta_{1}^2=\theta_{1}\in M$, $\theta_{5}^n \theta_{1}\in M$ for every $n\geq 0$. Since $\theta_{5}^n\theta_{4}\theta_{1}=\theta_{5}^{n+2}-\theta_{5}^n
			\theta_{1}\theta_{4}-\theta_{5}^{n+1}\theta_{1}\in M$ for every $n\geq 1$ and $\theta_{4}\theta_{1}\in M$ we have $\theta_{5}^n\theta_{4}\theta_{1}\in M$ for every $n\geq 0$.
			Since $\theta_{5}^n \theta_{4}\theta_{1}=\theta_{5}^{n+2}-\theta_{5}^n\theta_{1}\theta_{4}-\theta_{5}^{n+1}\theta_{1}$, for every $n\geq 1$, multiplying this equation by 
			$\theta_{1}$ on the right $\theta_{5}^n \theta_{4}\theta_{1}=\theta_{5}^{n+2}\theta_{1}-\theta_{5}^n\theta_{1}\theta_{4}\theta_{1}
			-\theta_{5}^{n+1}\theta_{1}$, for every $n\geq 1$. Then, $(\theta_{5}^n\theta_{1}\theta_{4})\theta_{1}=\theta_{5}^{n+2}\theta_{1}-\theta_{5}^n \theta_{4}\theta_{1}-\theta_{5}^{n+1}\theta_{1}\in M$. 	Since $(\theta_{5}^n\theta_{1})\theta_{1}=\theta_{5}^n\theta_{1}\in M$ for every $n\geq 1$,  $\theta_{3}\theta_{5}^n\theta_{1}=\theta_{1}\theta_{5}^n+\theta_{3}\theta_{5}^n
			-\theta_{3}\theta_{5}^{n-1}\theta_{4}-\theta_{5}^{n+1}\in M$ for every $n\geq 1$.
			
		 Furthermore, $(\theta_{1}\theta_{5}^n)\theta_{1}=\theta_{1}(\theta_{5}^n\theta_{1})\in \theta_{1}M\subset M$ for every $n\geq 1$ and $(\theta_{1}\theta_{5}^n\theta_{4})\theta_{1}=
		 \theta_{1}(\theta_{5}^n\theta_{4}\theta_{1})\in \theta_{1}M\subset M$ for every $n\geq 1$. Since $\theta_{3}\theta_{5}^n\theta_{4}\theta_{1}=-\theta_{1}\theta_{5}^n\theta_{4}-\theta_{1}\theta_{5}^{n+1}-\theta_{3}\theta_{5}^{n+1}+\theta_{5}^{n+1}\theta_{4}+
		 \theta_{3}\theta_{5}^{n+2}+\theta_{5}^{n+2}\in M$ for every $n\geq 0$ we have $M\theta_{1}\subset M$.
		 
		 \item $\theta_{5}M\subset M$.
		 
		 Note that $\theta_{5}(\theta_{4}\theta_{1})=(\theta_{5}\theta_{4})\theta_{1}\in M\theta_{1}\subset M$ since $\theta_{5}\theta_{4}\in M$.
		 On  the other hand, $\theta_{5}\theta_{3}=2\theta_{4}\theta_{1}-2\theta_{4}-2\theta_{5}^2+\theta_{5}\theta_{1}\in M$. Moreover, $\theta_{5}\theta_{1}\in M$, $\theta_{5}(\theta_{5}^n)=\theta_{5}^{n+1}\in M$ for every $n\geq 0$ and  $\theta_{5}(\theta_{5}^n\theta_{4})=\theta_{5}^{n+1}\theta_{4}\in M$ for every $n\geq 0$.
		 
		 However, $\theta_{5}(\theta_{5}^n\theta_{1}\theta_{4})=\theta_{5}^{n+1}\theta_{1}\theta_{4}\in M$ for every $n\geq 0$ and $\theta_{5}(\theta_{5}^n\theta_{1})
		 =\theta_{5}^{n+1}\theta_{1}\in M$ for $n\geq 1$. Furthermore, $\theta_{5}(\theta_{3}\theta_{5}^n)=(\theta_{5}\theta_{3})\theta_{5}^n\in M\theta_{5}^n\subset M$ for every $n\geq 1$, $\theta_{5}(\theta_{1}\theta_{5}^n)=(\theta_{5}\theta_{1})\theta_{5}^n \in M\theta_{5}^n\subset M$ for every $n\geq 1$, $\theta_{5}(\theta_{3}\theta_{5}^n\theta_{4})=
		 (\theta_{5}\theta_{3}\theta_{5}^n)\theta_{4}\in M\theta_{4}\subset M$ for every $n\geq 0$, $\theta_{5}(\theta_{1}\theta_{5}^n\theta_{4})=(\theta_{5}\theta_{1}\theta_{5}^n)\theta_{4}
		 \in M\theta_{4}\subset M$ for every $n\geq 1$. Thus, $\theta_{5}M\subset M$.
		 
		 \item  $\theta_{3} M\subset M$.
		 
		 Since $\theta_{3}\theta_{4}\in M$ we have that $\theta_{3}(\theta_{4}\theta_{1})=(\theta_{3}\theta_{4})\theta_{1}\in M\theta_{1}\subset M$. Since $\theta_{3}^2=
		 \theta_{3}-\theta_{5}+3\theta_{3}\theta_{4}+\theta_{3}\theta_{5}+\theta_{1}\theta_{4}+\theta_{5}\theta_{1}$ we have that $\theta_{3}^2\in M$ and $\theta_{3}\theta_{1}\in M$. Furthermore, 
		 $\theta_{3}\theta_{5}^n \in 	M$ for every $n\geq 0$, $\theta_{3}(\theta_{5}^n \theta_{4})=(\theta_{3}\theta_{5}^n)\theta_{4}\in M\theta_{4}\subset M$ for every $n\geq 0$.
		  Since $\theta_{3}(\theta_{5}^n\theta_{1}\theta_{4})=(\theta_{3}\theta_{5}^n)\theta_{1}\theta_{4}\in M\theta_{1}\theta_{4}\subset M\theta_{4}\subset M$ for every $n\geq 0$, 
		 $\theta_{3}(\theta_{5}^n\theta_{1})=(\theta_{3}\theta_{5}^n)\theta_{1}\in M\theta_{1}\subset M$ for every $n\geq 1$.
		 
		 On the other hand, $\theta_{3}(\theta_{3}\theta_{5})=\theta_{3}^2\theta_{5}^n\in M\theta_{5}^n\subset M$ for every $n\geq 1$ and $\theta_{3}(\theta_{1}\theta_{5}^n)=
		 (\theta_{3}\theta_{1})\theta_{5}^n\in M\theta_{5}^n\subset M$ for every $n\geq 1$. Since $\theta_{3}(\theta_{3}\theta_{5}^n \theta_{4})=\theta_{3}^2\theta_{5}^n\theta_{4}\in M\theta_{5}^n\theta_{4}\subset M\theta_{4} \subset M$ for every $n\geq 0$ and $\theta_{3}(\theta_{1}\theta_{5}^n\theta_{4})=(\theta_{3}\theta_{1})(\theta_{5}^n\theta_{4})\in 
		 M\theta_{5}^n\theta_{4}\subset M\theta_{4} \subset M$ for every $n\geq 1$ we have that $\theta_{3} M \subset M$.
		 
		 \item  $M\theta_{3} \subset M$.
		 
		 Note that $(\theta_{4}\theta_{1})\theta_{3}=\theta_{4}(\theta_{1}\theta_{3})\in \theta_{4}M\subset M$ and $\theta_{3}^2\in M$. Since $\theta_{3}\in  M$ we have that 
		 $\theta_{1}\theta_{3}\in \theta_{1}M\subset M$. On the other hand $\theta_{3} \in M$ implies $\theta_{5}^n \theta_{3}\in \theta_{5}^n M\subset M$ for every $n\geq 0$ and 
		 $(\theta_{5}^n \theta_{4})\theta_{3}=\theta_{5}^n (\theta_{4}\theta_{3})\in \theta_{5}^n M\subset M$ for every $n\geq 0$, since $\theta_{4}\theta_{3}\in M$. Note that 
		 $(\theta_{5}^n\theta_{1}\theta_{4})\theta_{3}=\theta_{5}^n\theta_{1}\theta_{4}\theta_{3}\in \theta_{5}^n\theta_{1}\theta_{4}M\subset \theta_{5}^n\theta_{1}M\subset
		 \theta_{5}^n M\subset M$ for every $n\geq 0$ and $(\theta_{5}^n \theta_{1})\theta_{3}=\theta_{5}^n \theta_{1}\theta_{3}\in \theta_{5}^n\theta_{1}M\subset
		 \theta_{5}^n M\subset M$ for every $n\geq 0$.	 Furthermore, since $\theta_{5}^n \theta_{3}\in M$ we have $(\theta_{3}\theta_{5}^n)\theta_{3}=\theta_{3}(\theta_{5}^n \theta_{3})\in \theta_{3}M\subset M$ for every $n\geq 1$ and 
		 $(\theta_{1}\theta_{5}^n)\theta_{3}=\theta_{1}(\theta_{5}^n \theta_{3}) \in \theta_{1}M \subset M$ for every $n\geq 1$.
		 
		 On the other hand, $\theta_{3}\theta_{5}^n\theta_{4}\theta_{3}=2\theta_{3}\theta_{5}^n\theta_{4}+\theta_{1}\theta_{5}^n\theta_{4}+\theta_{1}\theta_{5}^{n+1}
		 +\theta_{3}\theta_{5}^{n+1}+\theta_{3}\theta_{5}^{n+1}\theta_{4}-\theta_{5}^{n+1}\theta_{4}-\theta_{5}^{n+2}\in M$ for every $n\geq 0$ and 
		 $\theta_{1}\theta_{5}^n\theta_{4}=-\theta_{3}\theta_{5}^n\theta_{4}+2\theta_{5}^{n+1}\theta_{4}+\theta_{5}^{n+2}+\theta_{5}^{n+1}\theta_{1}\theta_{4}
		 +\theta_{3}\theta_{5}^{n+1}\theta_{4}+\theta_{5}^{n+2}\theta_{1}+\theta_{3}\theta_{5}^{n+2}+\theta_{1}\theta_{5}^{n+1}\theta_{4}-\theta_{1}\theta_{5}^{n+1}-
		 \theta_{3}\theta_{5}^{n+1} \in M$ for every $n\geq 0$. Thus,  $M\theta_{3} \subset M$. \qed
		 	\end{itemize}
\end{demos}

In \cite{VZ1} was proved that the algebra $\Gamma$ is generated by the elements \\
$$\alpha_{1}=\paren{\begin{matrix} 
	1 & 0 \\
	-1 & 0
	\end{matrix}}+\paren{\begin{matrix} 
	0 & 0 \\
	1 & 0
	\end{matrix}}z
+\paren{\begin{matrix} 
	1 & 1 \\
	0 & 0
	\end{matrix}}\frac{z^2}{2},$$ 
$$\alpha_{2}=\paren{\begin{matrix} 
	0 & 0 \\
	1 & 1
	\end{matrix}}+\paren{\begin{matrix} 
	0 & 0 \\
	-1 & 0
	\end{matrix}}z
+\paren{\begin{matrix} 
	-1 & -1 \\
	0 & 0
	\end{matrix}}\frac{z^2}{2},$$
$$\alpha_{3}=\paren{\begin{matrix} 
	1 & 1 \\
	-1 & -1
	\end{matrix}}z
+\paren{\begin{matrix} 
	-1 & 1 \\
	0 & 0
	\end{matrix}}\frac{z^2}{2}, \hspace{0.1cm}
	\alpha_{4}=\paren{\begin{matrix} 
	0 & 0 \\
	1 & 0
	\end{matrix}}\frac{z^2}{2}, \hspace{0.1cm}
		\text{and} \hspace{0.1cm} \alpha_{5}=\paren{\begin{matrix} 
	0 & 0 \\
	0 & 1
	\end{matrix}}\frac{z^2}{2}.$$

In the following proposition we look for the generators which give us the presentation. 
\begin{prop}\label{caloind}
Define $\beta_{1}=\alpha_{1}+\alpha_{3}$, $\beta_{3}=\alpha_{1}-\alpha_{3}$, $\beta_{4}=2\alpha_{4}$, $\beta_{5}=2\alpha_{5}$ then 
	$\llave{\beta_{4}\beta_{1},\beta_{3},\beta_{1}}\cup$ \\ $\llave{\beta_{5}^n\mid n\geq 0}\cup \llave{\beta_{5}^n \beta_{4} \mid n\geq 0}\cup
	\llave{\beta_{5}^n \beta_{1}\beta_{4}  \mid n\geq 0}\cup \llave{\beta_{5}^n \beta_{1} \mid n\geq 1}\cup
	\llave{\beta_{3}\beta_{5}^n \mid n\geq 1}\cup \llave{\beta_{1}\beta_{5}^n \mid  n\geq 1}\cup$  \\ $\llave{\beta_{3}\beta_{5}^n\beta_{4} \mid n\geq 0}
	\cup \llave{\beta_{1}\beta_{5}^n\beta_{4} \mid n\geq 1}$ is a linearly independent set over $\mathbb{K}$.
\end{prop}
\begin{demos}
	Note that 
	
	\begin{equation*}
	\mathbb{A}=\mathbb{K} \cdot \beta_{1} \oplus \mathbb{K} \cdot \beta_{3} \oplus \mathbb{K}\cdot\beta_{4}\beta_{1} \oplus 
	\bigoplus_{n=0}^{\infty} \mathbb{K} \cdot \beta_{5}^n\oplus
	\bigoplus_{n=0}^{\infty} \mathbb{K} \cdot\beta_{5}^n \beta_{4} \oplus 
	\bigoplus_{n=0}^{\infty} \mathbb{K}\cdot \beta_{5}^n \beta_{1}\beta_{4} \oplus 
	\bigoplus_{n=1}^{\infty} \mathbb{K}\cdot \beta_{5}^n \beta_{1} \oplus 
	\bigoplus_{n=1}^{\infty} \mathbb{K} \cdot \beta_{3}\beta_{5}^n \oplus 
		\end{equation*}
			\begin{equation*}
	\bigoplus_{n=1}^{\infty} \mathbb{K} \cdot \beta_{1}\beta_{5}^n \oplus 
	\bigoplus_{n=1}^{\infty} \mathbb{K} \cdot \beta_{3}\beta_{5}^n\beta_{4}\oplus 
	\bigoplus_{n=1}^{\infty} \mathbb{K} \cdot \beta_{1}\beta_{5}^n\beta_{4} =
	\mathbb{K}\oplus	\mathbb{K} \cdot \beta_{1} \oplus \mathbb{K} \cdot \beta_{3} 
		\end{equation*}
			\begin{equation*}
	\oplus \mathbb{K}\cdot\paren{\frac{1}{2}\beta_{4}\beta_{1} -\frac{1}{2}\beta_{4}\beta_{3}+\beta_{4}
		+\beta_{5}\beta_{1}+\frac{1}{2}\beta_{5}\beta_{4}-\frac{1}{2}\beta_{5}^2-\beta_{3}\beta_{4}}
	\end{equation*}
	\begin{equation*}
	\oplus 
	\mathbb{K}\cdot \paren{	\frac{1}{2}\beta_{4}\beta_{1} -\frac{1}{2}\beta_{4}\beta_{3}+\beta_{4}+\beta_{5}\beta_{1}+\frac{1}{2}\beta_{5}\beta_{4}-\frac{1}{2}\beta_{5}^2 }	
	\oplus 
	\bigoplus_{n=0}^{\infty} \mathbb{K} \cdot (\beta_{5}^n-\beta_{5}^{n-1} \beta_{4} -\beta_{3}\beta_{5}^n  )\oplus
	\bigoplus_{n=0}^{\infty} \mathbb{K} \cdot \beta_{5}^n\beta_{4} 
				\end{equation*}
			\begin{equation*}
	\oplus 
	\bigoplus_{n=0}^{\infty} \mathbb{K} \cdot (\beta_{5}^n \beta_{1}\beta_{4} +\beta_{3}\beta_{5}^n\beta_{4})\oplus 
	\bigoplus_{n=1}^{\infty} \mathbb{K} \cdot \beta_{5}^n \beta_{1}+\beta_{5}^{n-1} \beta_{4} \oplus 
	\bigoplus_{n=1}^{\infty} \mathbb{K} \cdot \beta_{5}^n \beta_{1}\beta_{4} 
		\end{equation*}
	\begin{equation*}
	 \oplus 
	\bigoplus_{n=1}^{\infty} \mathbb{K} \cdot (\beta_{1}\beta_{5}^n + \beta_{3}\beta_{5}^n )\oplus 
	\bigoplus_{n=1}^{\infty} \mathbb{K} \cdot \beta_{5}^n \oplus 
	\bigoplus_{n=1}^{\infty} \mathbb{K} \cdot (\beta_{1}\beta_{5}^n\beta_{4} +\beta_{3}\beta_{5}^n\beta_{4}).
	\end{equation*}

	The second equality is given by an isomorphism of $\mathbb{K}$ vector spaces  sending 
	
	$\llave{\beta_{4}\beta_{1},\beta_{3},\beta_{1}}\cup \llave{\beta_{5}^n\mid n\geq 0}
	\cup \llave{\beta_{5}^n \beta_{4} \mid n\geq 0}\cup
	\llave{\beta_{5}^n \beta_{1}\beta_{4}  \mid n\geq 0}\cup \llave{\beta_{5}^n \beta_{1} \mid n\geq 1}$
	
	$\cup
	\llave{\beta_{3}\beta_{5}^n \mid n\geq 1}\cup \llave{\beta_{1}\beta_{5}^n \mid  n\geq 1}\cup \llave{\beta_{3}\beta_{5}^n\beta_{4} \mid n\geq 0}
	\cup \llave{\beta_{1}\beta_{5}^n\beta_{4} \mid n\geq 1}$ to the set $$\llave{1,\beta_{1},\beta_{3},\frac{1}{2}\beta_{4}\beta_{1} -\frac{1}{2}\beta_{4}\beta_{3}+\beta_{4}
		+\beta_{5}\beta_{1}+\frac{1}{2}\beta_{5}\beta_{4}-\frac{1}{2}\beta_{5}^2-\beta_{3}\beta_{4},
		\frac{1}{2}\beta_{4}\beta_{1} -\frac{1}{2}\beta_{4}\beta_{3}+\beta_{4}+\beta_{5}\beta_{1}
		+\frac{1}{2}\beta_{5}\beta_{4}-\frac{1}{2}\beta_{5}^2 }$$
$$	\cup \llave{ \beta_{5}^n-\beta_{5}^{n-1} \beta_{4} -\beta_{3}\beta_{5}^n \mid
		n\geq 0} \cup \llave{ \beta_{5}^n\beta_{4} \mid n\geq 0  } \cup \llave{\beta_{5}^n \beta_{1}\beta_{4} +\beta_{3}\beta_{5}^n\beta_{4}\mid n\geq 0  }\cup
	\llave{\beta_{5}^n \beta_{1}+\beta_{5}^{n-1} \beta_{4} \mid n\geq 1}$$
	$$\cup \llave{ \beta_{5}^n \beta_{1}\beta_{4} \mid n\geq 1} \cup
	\llave{\beta_{1}\beta_{5}^n + \beta_{3}\beta_{5}^n\mid n\geq 1 }\cup \llave{\beta_{5}^n \mid n\geq 1} \cup \llave{\beta_{1}\beta_{5}^n\beta_{4} +\beta_{3}\beta_{5}^n\beta_{4}\mid
		n \geq 1}$$ which is linearly independent because is exactly $\llave{1,\beta_{1},\beta_{3},\beta_{4},\beta_{5}}\cup \llave{e_{ij}x^k \mid 1\leq i,j \leq 2, k\geq 3}$. \qed

	Finally, we conclude with the proof of the  presentation.  
	Define 
	\begin{equation*}
	f:\complex\cdot \langle \theta_{1},\theta_{3}, \theta_{4}, \theta_{5}\rangle /I \longrightarrow \mathbb{A},
	\end{equation*}
	\begin{equation*}
	f(\overline{\theta_{j}})=\beta_{j}.
	\end{equation*}
The Lemma \ref{calogen} guarantees the existence of  the system of generators
	$$\llave{\theta_{4}\theta_{1},\theta_{3},\theta_{1}}\cup \llave{\theta_{5}^n\mid n\geq 0}\cup \llave{\theta_{5}^n \theta_{4} \mid n\geq 0}\cup
	\llave{\theta_{5}^n \theta_{1}\theta_{4}  \mid n\geq 0}\cup \llave{\theta_{5}^n \theta_{1} \mid n\geq 1}\cup
	\llave{\theta_{3}\theta_{5}^n \mid n\geq 1}\cup$$ $$\llave{\theta_{1}\theta_{5}^n \mid  n\geq 1}\cup \llave{\theta_{3}\theta_{5}^n\theta_{4} \mid n\geq 0}
	\cup \llave{\theta_{1}\theta_{5}^n\theta_{4} \mid n\geq 1}$$   for $\complex\cdot \langle \theta_{1},\theta_{3}, \theta_{4}, \theta_{5}\rangle /I $ as a free  $\complex$-vector space. Furthermore
	$f\mid _{\complex}: \complex\longrightarrow A$ is a monomorphism.
	
	The Proposition \ref{caloind} implies that 	$\llave{\beta_{4}\beta_{1},\beta_{3},\beta_{1}}\cup \llave{\beta_{5}^n\mid n\geq 0}\cup \llave{\beta_{5}^n \beta_{4} \mid n\geq 0}\cup
	\llave{\beta_{5}^n \beta_{1}\beta_{4}  \mid n\geq 0}\cup \llave{\beta_{5}^n \beta_{1} \mid n\geq 1}\cup
	\llave{\beta_{3}\beta_{5}^n \mid n\geq 1}\cup \llave{\beta_{1}\beta_{5}^n \mid  n\geq 1}\cup \llave{\beta_{3}\beta_{5}^n\beta_{4} \mid n\geq 0}
	\cup \llave{\beta_{1}\beta_{5}^n\beta_{4} \mid n\geq 1}$ is a linearly independent set over $\complex$. 
	
	Putting together Lemma~\ref{calogen} and  Propositions~\ref{caloind} we conclude the proof of Theorem~\ref{calogero}. \qed
\end{demos}
\end{demos}

\section{Conclusion and Final Comments}

In this article, we obtained in  Theorem~\ref{presentation} a general result for presentations of finitely generated algebras. The theorem can be used to obtain a complete description in terms of generators and relations since it says when a set of relations is enough to characterize a given finitely generated algebra. As an application, we find nice presentations for matrix bispectral algebras and give positive answers for the conjectures presented in~\cite{Grfrm-4}. 

An important role was played by the Ad-condition due to the fact that the matrix-valued operators were acting from opposite directions, since we can consider these algebras as matrix polynomial.


Another research direction would be to investigate the presentations of the full rank $1$ algebras which by Theorem 1 in \cite{VZ1} are finitely generated. As we saw, the examples given in \cite{Grfrm-4} and worked out here, are finitely presented. However, this is not necessarily true for general non-commutative rings.

\section*{Acknowledgments}

BVDC was supported by CAPES grants 88882 332418/2019-01 as well as IMPA.
JPZ was supported by CNPq grants 302161 and 47408, as well as by FAPERJ under the program \textit{Cientistas do Nosso Estado} grant E-26/202.927/2017, Brazil. BDVC and JPZ acknowledge the support from the FSU-2020-09 grant from Khalifa University, UAE.


\begin{thebibliography}{10}

\bibitem{duistermaat1986differential}
J.J. Duistermaat and F.A. Gr{\"{u}}nbaum.
\newblock Differential equations in the spectral parameter.
\newblock {\em Communications in Mathematical Physics}, 103(2):177--240, 1986.

\bibitem{Zubelli1991}
Jorge~P. Zubelli and Franco Magri.
\newblock Differential equations in the spectral parameter, {D}arboux
  transformations and a hierarchy of master symmetries for {K}d{V}.
\newblock {\em Comm. Math. Phys.}, 141(2):329--351, 1991.

\bibitem{marta2018}
Tom~H. {Koornwinder} and Marta {Mazzocco}.
\newblock {Dualities in the \(q\)-Askey scheme and degenerate DAHA}.
\newblock {\em {Stud. Appl. Math.}}, 141(4):424--473, 2018.

\bibitem{AMM1977}
H.~{Airault}, H.~P. {McKean}, and J\"urgen {Moser}.
\newblock {Rational and elliptic solutions of the Korteweg-de Vries equation
  and a related many body problem}.
\newblock {\em {Commun. Pure Appl. Math.}}, 30:95--148, 1977.

\bibitem{wilson}
George Wilson.
\newblock Bispectral commutative ordinary differential operators.
\newblock {\em J. Reine Angew. Math.}, 442:177--204, 1993.

\bibitem{Chalub2000}
Fabio A. C.~C. Chalub and Jorge~P. Zubelli.
\newblock Integrable systems, {H}uygens' principle, and {D}irac operators.
\newblock In {\em Proceedings of the Workshop on Nonlinearity, Integrability
  and All That: Twenty Years after NEEDS '79 (Gallipoli, 1999)}, pages 89--96,
  River Edge, NJ, 2000. World Sci. Publishing.

\bibitem{Chalub2001}
Fabio A. C.~C. Chalub and Jorge~P. Zubelli.
\newblock On {H}uygens' principle for {D}irac operators and nonlinear evolution
  equations.
\newblock {\em J. Nonlinear Math. Phys.}, 8(suppl.):62--68, 2001.
\newblock Nonlinear evolution equations and dynamical systems (Kolimbary,
  1999).

\bibitem{Chalub2001a}
F.~A. C.~C. Chalub and J.~P. Zubelli.
\newblock S{\'o}litons: Na crista da onda por mais de 100 anos.
\newblock {\em Revista Matem{\'a}tica Universit{\'a}ria}, 30:44--62, 2001.

\bibitem{MR2201201}
Fabio A. C.~C. Chalub and Jorge~P. Zubelli.
\newblock Huygens' principle for hyperbolic operators and integrable
  hierarchies.
\newblock {\em Phys. D}, 213(2):231--245, 2006.

\bibitem{Sakhnovich2001}
A.~L. Sakhnovich and J.~P. Zubelli.
\newblock Bundle bispectrality for matrix differential equations.
\newblock {\em Integral Equations and Operator Theory}, 41(4):472--496, 2001.

\bibitem{PNAS2019}
W.~Riley Casper, F.~Alberto Gr{\"u}nbaum, Milen Yakimov, and Ignacio
  Zurri{\'a}n.
\newblock Reflective prolate-spheroidal operators and the kp/kdv equations.
\newblock {\em Proceedings of the National Academy of Sciences},
  116(37):18310--18315, 2019.

\bibitem{fokas1987}
A.~S. {Fokas}.
\newblock {Symmetries and integrability}.
\newblock {\em {Stud. Appl. Math.}}, 77:253--299, 1987.

\bibitem{fokas2002}
F.~{Finkel} and A.~S. {Fokas}.
\newblock {On the construction of evolution equations admitting a master
  symmetry}.
\newblock {\em {Phys. Lett., A}}, 293(1-2):36--44, 2002.

\bibitem{DERKSEN2015210}
Harm Derksen and Jiarui Fei.
\newblock General presentations of algebras.
\newblock {\em Advances in Mathematics}, 278:210--237, 2015.

\bibitem{zbMATH03512288}
Ren\'ee {Elkik}.
\newblock {Solutions d'\'equations \`a coefficients dans un anneau
  hens\'elien}.
\newblock {\em {Ann. Sci. \'Ec. Norm. Sup\'er. (4)}}, 6:553--603, 1973.

\bibitem{zbMATH01713734}
Alberto {Arabia}.
\newblock {Rel\`evements des alg\`ebres lisses et de leurs morphismes}.
\newblock {\em {Comment. Math. Helv.}}, 76(4):607--639, 2001.

\bibitem{Grfrm-4}
F.~Alberto Gr{\"{u}}nbaum.
\newblock Some noncommutative matrix algebras arising in the bispectral
  problem.
\newblock {\em SIGMA Symmetry Integrability Geom. Methods Appl.}, 10:078, 2014.

\bibitem{zbMATH05249007}
Mirta~M. {Castro} and F.~Alberto {Gr\"unbaum}.
\newblock {The algebra of differential operators associated to a family of
  matrix-valued orthogonal polynomials: five instructive examples}.
\newblock {\em {Int. Math. Res. Not.}}, 2006(7):33, 2006.
\newblock Id/No 47602.

\bibitem{zbMATH01684157}
F.~Alberto Gr{\"{u}}nbaum.
\newblock {The bispectral problem: An overview}.
\newblock In {\em Special functions 2000: current perspective and future
  directions. Proceedings of the NATO Advanced Study Institute, Tempe, AZ, USA,
  May 29--June 9, 2000}, pages 129--140. Dordrecht: Kluwer Academic Publishers,
  2001.

\bibitem{tirao}
Tirao J.
\newblock The algebra of differential operators associated to a weight matrix:
  a first example, in groups, algebras and applications.
\newblock {\em Contemporary Mathematics}, 537:291--324, 2011.

\bibitem{VZ1}
Brian Vasquez and Jorge Zubelli.
\newblock Matrix bispectrality of full rank one algebras.
\newblock {\em axXiv preprint}, 2021.

\bibitem{BGK09}
Maarten Bergvelt, Michael Gekhtman, and Alex Kasman.
\newblock Spin calogero particles and bispectral solutions of the matrix {KP}
  hierarchy.
\newblock {\em Mathematical Physics Analysis and Geometry}, 12, 07 2008.

\bibitem{zbMATH06722531}
Joel {Geiger}, Emil {Horozov}, and Milen {Yakimov}.
\newblock {Noncommutative bispectral Darboux transformations}.
\newblock {\em {Trans. Am. Math. Soc.}}, 369(8):5889--5919, 2017.

\bibitem{zurrian}
Ignacio Zurri\'{a}n.
\newblock {The Algebra of Differential Operators for a Matrix Weight: An
  Ultraspherical Example}.
\newblock {\em International Mathematics Research Notices}, 2017(8):2402--2430,
  06 2016.

\end{thebibliography}

\end{document}